\documentstyle[a4,11pt,french]{article}

\begin{document}

\newtheorem{prop}[subsection]{Proposition}
\newtheorem{lemm}[subsection]{Lemme}
\newtheorem{cor}[subsection]{Corollaire}
\newtheorem{theo}[subsection]{Th\'{e}or\`{e}me}
\newcounter{cms}

\def\grdZ{{\bf Z}}
\def\grdF{{\bf F}}
\def\grdT{{\bf T}}
\def\grdR{{\bf R}}
\def\grdC{{\bf C}}
\def\grdN{{\bf N}}
\def\grdQ{{\bf Q}}
\def\grdP{{\bf P}}
\def\grdG{{\bf G}}
\def\grdH{{\bf H}}

\title{\bf{Bornes effectives pour la torsion des courbes elliptiques sur
les corps de nombres.}}
\author{Pierre Parent}
\date{8 novembre 1996}

\maketitle

\begin{abstract}
On se propose de donner une forme effective au th\'{e}or\`{e}me de
Mazur-Kamienny-Merel sur la torsion des courbes elliptiques sur les corps de
nombres.
\end{abstract}
\tableofcontents
\section{Pr\'esentation des r\'esultats.}
\subsection{Introduction.}
{\it  La ``conjecture de borne uniforme pour les courbes
elliptiques'', affirmant qu'il existe pour tout entier $d$ un entier $B(d)$
tel que, pour tout corps de nombres $K$ de degr\'{e} $d$ sur $\grdQ $ et pour
toute courbe elliptique $E$ sur $K$, la partie de torsion $E(K)_{\rm tors}$
du groupe de Mordell-Weil
$E(K)$ est de cardinal major\'{e} par $B(d)$, a \'{e}t\'{e} d\'{e}montr\'{e}e
dans le cas g\'{e}n\'{e}ral en f\'{e}vrier 1994 par Lo\"{\i}c Merel. En
fait, Merel (et Oesterl\'{e})
montrent que, si $P$ est un point d'ordre $p$ premier de $E(K)$, on a
$p \leq (1+{3^{d/2}})^2 $ . Des travaux de Faltings et Frey permettent alors de
conclure \`{a} l'existence des bornes $B(d)$, mais pas de mani\`{e}re
effective : en effet, si on a bien major\'{e} les nombres premiers pouvant
diviser les groupes $E(K)_{\rm tors}$, on ne sait pas en pratique quelles
puissances de ces nombres premiers peuvent intervenir dans ces groupes.
  Le but de cet article est de d\'emontrer une forme explicite de la forme
forte de la conjecture de borne uniforme (le corollaire 1.8 ci-dessous), en
donnant une borne pour ces puissances de premiers qui peuvent diviser la
torsion.}
\subsection{Sch\'{e}ma de la preuve.}
Soit $E$ une courbe elliptique sur un corps de nombres $K$
de degr\'{e} $d$ sur $\grdQ $, poss\'{e}dant un point $K$-rationnel $P$, tel
que le cardinal du sous-groupe cyclique de $E(K)$ engendr\'{e} par $P$ soit
une puissance d'un nombre premier : $|\langle P\rangle |={p^n}$. On cherche
\`{a} majorer $p^n $ en fonction de $d$. Soit $l$ un nombre premier
diff\'{e}rent de $2$ et de $p$ (on prendra dans la suite le plus petit
possible, $i.e.\ l=5$ si $p=3$, $l=3$ dans tous les autres cas). Soit $\cal L$
un id\'eal maximal de ${\cal O}_K$ (l'anneau des entiers du corps $K$),
contenant $l$. Examinons la fibre en $\cal L$ du mod\`{e}le de N\'{e}ron
(qu'on note ${\cal E}$) de la courbe elliptique sur ${\cal O}_K$.
\begin{prop}
Si en $\cal L$ au-dessus de $l$, on a l'un des trois cas :
\begin{enumerate}
\item ${\cal E}$ a bonne r\'eduction ;
\item ${\cal E}$ a r\'eduction additive ;
\item ${\cal E}$ a r\'eduction multiplicative tordue ;
\item ${\cal E}$ a r\'eduction multiplicative d\'eploy\'ee et $\langle
\widetilde{P} \rangle$ (o\`u $\widetilde{P}$ est la r\'{e}duction de $P \bmod
\cal L$) appartient \`a la composante neutre,
\end{enumerate}
alors $|\langle P\rangle |\leq 2.(1+l^d )$.
\end{prop}
{\bf Remarque.} {\it  Cette proposition classique est explicit\'ee dans
\cite{Bas} par exemple, ou dans \cite{kamienny}. Dans le cas 2., on peut
borner %
par $4$ l'ordre du groupe des composantes de la r\'eduction additive, donc la
borne pour l'ordre de $P$ est $1$ si $p\geq 5$, $3$ si $p=3$, et $4$ si $p=2$.
Dans les cas 1. et 4., la borne est en fait $(l^{d/2} +1)^2$ (borne de Weil)
et $(l^d -1)$ (cardinal de groupe multiplicatif) respectivement, donc en
$l^d$. Le cas 3. la porte \`a $(1+l^d )$ ou \`a $2(1+3^d )$ : dans ce cas en
effet, soit $\widetilde{P}$ appartient \`a la composante neutre, et son ordre
est major\'e par $(1+l^d)$ ; soit il est dans une composante non triviale. Le
groupe de Galois d'une extension quadratique du corps r\'esiduel $k({\cal L})$
agit alors par la multiplication par $(\pm 1)$ sur le groupe des composantes.
Mais $\widetilde{P}$ est $k({\cal L})$-rationnel, et donc la composante \`a
laquelle il appartient \'egale son oppos\'ee, ce qui veut dire que
$2\widetilde{P}$ est dans la composante neutre (chose qui ne peut arriver que
si $p$ est $2$).}

   Supposons donc qu'en tout $\cal L$ au-dessus de $l$, ${\cal
E}$ ait r\'eduction multiplicative d\'eploy\'ee et que $\widetilde{P}$ ne soit
pas trivial dans le groupe des composantes de cette r\'eduction ${\cal
E}_{k({\cal L})}$ de ${\cal E}_{/{\cal O}_K}$. Pour avoir une bonne
interpr\'etation modulaire dans la suite, comme indiqu\'e plus bas, on
aimerait que $\widetilde{P}$ soit d'ordre $p^n$ dans le groupe des composantes
de ${\cal E}_{k({\cal L})}$ : d'o\`u l'id\'ee d'examiner le quotient de
${\cal E}_{k({\cal L})}$ par le plus gros sous-groupe de $\langle
\widetilde{P} \rangle$ inclus dans sa composante neutre ${\cal E}^0_{k({\cal
L})}$. Soit donc $n_{\cal L}$ le plus petit entier tel que ${p^{n_{\cal L}} }.
\widetilde{P}$ tombe dans ${\cal E}^0_{k({\cal L})}$. Soit aussi $n'$ le plus
petit des $n_{\cal L}$, pour ${\cal L}$ parcourant l'ensemble des places de
${\cal O}_K$ au-dessus de $l$. Oesterl\'{e} \'enonce :
\begin{lemm}
Soit $k$ \'{e}l\'{e}ment de $\grdN $ ; si en la place $\cal L$, $p^{k-1} .P$
ne se r\'{e}duit pas dans la composante neutre de ${\cal E}_{k({\cal L})}$,
alors la r\'{e}duction de $\widetilde{P}$ dans la fibre en ${\cal L}$ du
mod\`ele de N\'eron ${\cal E}'$ de $E/\langle {p^k}.P\rangle$ est d'ordre
exactement $p^k$ dans son groupe des composantes.
\end{lemm}
L'interpr\'{e}tation modulaire de ce lemme, qui sera d\'emontr\'e dans la
section~2, est donc la suivante : le point $K-$rationnel $j$  de la courbe
modulaire $X_0 (p^{n'})$ que d\'efinit le couple $(E/\langle {p^{n'}}.P\rangle
,\langle P\rangle )$ se r\'{e}duit en la pointe $0$ modulo toute place
$\cal L$ ; et l'image $j'$ de ce point par l'involution d'Atkin-Lehner, en la
pointe infinie (voir \cite{mazur} , page 159).  Voyons comment on peut alors
utiliser les arguments de Kamienny (\cite{kamienny}, \cite{Bas}) pour
d\'emontrer son crit\`{e}re.\\
{\bf Remarque.} {\it Dans toute la suite, on notera $E$ ce qui sera en
r\'ealit\'e $(E/\langle {p^{n'}}.P\rangle )$ ; on notera aussi $p^n$ pour
$p^{n'}$. La borne qu'on obtiendra pour $|\langle P\rangle |$ devra donc
\^etre multipli\'ee par un facteur $(l^d -1)$ pour avoir une borne \`a notre
$p^n$ ``originel'' (en effet, $(p^n /p^{n'})$ est inf\'erieur au cardinal du
groupe multiplicatif ${\grdG}_{m,{\grdF}_l}$)}.

  Si les ${\sigma }_i,\ 1\leq i\leq d$ sont les plongements de $K$ dans
$\overline{\grdQ }$, ${j'}^{(d)} :=({\sigma}_1 (j'),{\sigma }_2 (j'),\\
..., {\sigma}_d (j'))$ d\'{e}finit un point $\grdQ $-rationnel du produit
sym\'{e}trique $d$-i\`{e}me : ${X_0 (p^n)}^{(d)}$, de $X_0 (p^n)$.

  D\'efinissons comme Merel le quotient d'enroulement (\cite{merel}).
On consi\-d\`{e}re les premiers groupes d'homologie singuli\`{e}re absolue  :
$H_1 (X_0 (p^n )\, ;\, \grdZ )$ et relative aux pointes :
$H_1 (X_0 (p^n ),\, {\rm pointes}\, ;\, \grdZ )$, de $X_0(p^n )$, le
premier \'etant vu comme un sous-groupe du second. Si $a$ et $b$ sont deux
\'el\'ements de $\grdP^1(\grdQ)$, le {\it symbole modulaire} $\{a,b\}$ est
l'\'{e}l\'{e}ment de $H_1 (X_0 (p^n ),\, {\rm pointes}\, ;\, \grdZ )$
d\'efini par l'image de n'importe quel chemin continu reliant $a$ \`{a} $b$
sur le demi-plan de Poincar\'{e} auquel on a ajout\'{e} l'ensemble
${\grdP }^1 (\grdQ )$ de ses pointes. L'int\'egration d\'efinit un
isomorphisme classique d'espaces vectoriels r\'eels : %
$$\left\{ \begin{array}{c}
H_1 (X_0 (p^n ) \, ;\, \grdZ )\otimes \grdR \to
{\rm Hom}_{\grdC } \left( H^0 (X_0 (p^n )\, ;\, {\Omega }^1 ),
\ \grdC \right) \\
\gamma \otimes 1 \mapsto \left( \omega \mapsto \int_{\gamma } \omega \right).
\end{array} \right. $$
Selon un th\'eor\`eme de Manin et Drinfeld, %
l'image %
r\'eciproque de la forme lin\'{e}aire $\omega \mapsto \int_{\{ 0,\infty \} }
\omega $ dans $H_1 (X_0 (p^n )\, ;\, \grdR )$ est en r\'{e}alit\'{e}  dans
$H_1 (X_0 (p^n )\, ;\, \grdQ )$. C'est {\it l'\'el\'ement d'enroulement},
qu'on note $e$ (comme d'habitude). (En fait, le r\'esultat de Manin assure
plus g\'en\'eralement que l'image r\'eciproque de l'int\'egration sur {\em
tout} symbole modulaire est \`a coefficients dans $\grdQ$.) Notons (toujours
comme d'habitude) $\grdT$ l'alg\`ebre engendr\'ee sur $\grdZ$ par les
op\'erateurs de Hecke $T_i$ ($i\geq 1$, entier), agissant fid\`element entre
autres sur $H_1 (X_0 (p^n )\, ;\, \grdQ )$ et sur la jacobienne $J_0 (p^n )$
de la courbe modulaire. Soit ${\cal A}_e$ l'id\'eal annulateur dans $\grdT$ de
$e$ ({\it id\'eal d'enroulement}) ; on d\'{e}finit alors le {\em quotient
d'enroulement} $J_0^e$ comme la vari\'et\'e ab\'elienne quotient $J_0 (p^n )/
{{\cal A}_e} J_0 (p^n )$. De fa\c{c}on similaire \`a Merel dans \cite{merel},
un th\'eor\`eme de Kolyvagin-Logachev nous permet de montrer dans la section~3
le : %
\begin{theo}
$J_0^e (\grdQ )$ est fini.
\end{theo}
Soit maintenant l'application naturelle $f_d$ :
${X_0 (p^n )}^{(d)}_{\rm lisse} \rightarrow J_0^e $,
qu'on a normalis\'ee par ${\infty }^{(d)}
\mapsto 0$. On montre ``par les arguments standards'' que le fait que $f_d$
soit une immersion formelle en $\infty _{{\grdF}_l}^{(d)}$ contredit
l'existence de notre point $j'^{(d)}$ (voir la sous-section 4.12). Or on a le
``crit\`ere de Kamienny'' :
\begin{theo}
 On a \'{e}quivalence entre :
\begin{enumerate}
\item $f_d$ est une immersion formelle en ${\infty}_{{\grdF}_l}^{(d)}$, et
\item $T_1 e,...,T_d e$ sont $\grdF _l$-lin\'{e}airement ind\'{e}pendants dans
$\grdT e/l\grdT e$.
\end{enumerate}
De plus, ces deux conditions sont satisfaites si l'est :
\begin{enumerate}
\item[3.] $T_1 \{ 0,\infty \} ,...,T_{d.s} \{ 0,\infty \}$ sont
$\grdF_l$-lin\'{e}airement ind\'{e}pendants dans l'espace vectoriel
$H_1 (X_0 (p^n ),\, {\mathrm {pointes }}\, ;\, \grdZ )\otimes \grdF_l$ (ici et
pour la suite, $s$ d\'{e}signe le plus petit nombre premier diff\'erent de
$p$).
\end{enumerate}
\end{theo}
(Le fait que la derni\`{e}re condition implique les pr\'{e}c\'{e}dentes est
une remarque d'Oesterl\'{e} en niveau premier, utilis\'ee d\'ej\`a par Merel ;
ce th\'eor\`eme sera d\'emontr\'e dans la section~4). Il suffit donc
maintenant de prouver :
\begin{prop}
Soit $C:=\sqrt{65}$, si $p$ est diff\'erent de $2$, et $C:=\sqrt{129}$ si $p$
est $2$. Si $p^n\geq C^2 .(sd)^6$, alors les $T_i \{ 0,\infty \} ,\ 1\leq i
\leq sd$ sont $\grdF$-lin\'{e}airement ind\'{e}pendants (dans le
$\grdF$-espace vectoriel $H_1 (X_0 (p^n ),\, {\mathrm {pointes }}\, ;\, \grdZ
)\otimes \grdF$) pour tout corps $\grdF$.
\end{prop}
De cette proposition, dont la d\'emonstration occupe la section~5, d\'ecoule
donc le :
\begin{cor}
Soit $E$ une courbe elliptique sur un corps  $K$ de degr\'{e} $d$ sur
$\grdQ$. Si $E(K)$ poss\`{e}de un point $P$ d'ordre une puissance $p^n$ d'un
nombre premier $p$, on a :
\begin{enumerate}
\item $p^n \leq 65. (3^d -1).(2d)^6 \ ,$ si $p$ est diff\'{e}rent de $2$
et $3$ ;
\item Si $p=3$, $p^n \leq 65.(5^d -1).(2d)^6 $ ;
\item et pour $p=2$, $2^n \leq 129.(3^d -1).(3d)^6 .$
\end{enumerate}
\end{cor}
\noindent{\it {\bf Remarques.} Le th\'eor\`eme de Mordell-Weil assure que pour
tous $K$ et $E$ comme ci-dessus, il existe deux entiers $n_1$ et $n_2$ tel que
$E(K)_{\rm tors} \cong  {\grdZ }/{{n_1}\grdZ }\times {\grdZ }/{{n_2}
\grdZ }$. On borne donc ainsi le cardinal de tous les $p$-groupes de
$E(K)_{\rm tors}$, et avec la borne de Merel-Oesterl\'e pour les nombres
premiers $p$ pouvant intervenir on obtient une borne effective ``globale''
pour l'ensemble $\{ card( E(K)_{\rm tors})\  |\  K$ est un corps de
nombres de degr\'e $d$ et $E$ est une courbe elliptique sur $K\}$.

  Une telle borne globale semble de toute fa\c{c}on assez facile \`{a}
am\'{e}liorer en reprenant les arguments de ce papier directement
en niveau $N$ entier quelconque.}

  La preuve dans son ensemble ayant \'{e}t\'{e} esquiss\'{e}e, nous allons
montrer dans la suite, dans l'ordre indiqu\'e, les diff\'erentes propositions
utilis\'ees.
\section{Quotients de courbes elliptiques.}
On prouve le lemme 1.4 :
\begin{prop}
Soit $E$ une courbe elliptique sur un corps de nombres $K$, notons ${\cal E}$
son mod\`ele de N\'eron sur ${\cal O}_K$, et soit $P$ un point $K$-rationnel
de $E$, d'ordre une puissance $p^n$ d'un nombre premier $p$. Soit encore $k$
un entier, et ${\cal L}$ une place de ${\cal O}_K$ qui ne soit pas au-dessus
de $p$. Supposons qu'en ${\cal L}$, ${\cal E}$ ait r\'eduction multiplicative
d\'eploy\'ee, et que $p^{k-1} .P$ ne se r\'eduise pas dans la composante
neutre de ${\cal E}_{k({\cal L})}$. Alors la r\'eduction de $P$ dans la fibre
en ${\cal L}$ du mod\`ele de N\'eron ${\cal E}'$ de $E/\langle p^k .P\rangle$
est d'ordre exactement $p^k$ dans son groupe de composantes.
\end{prop}
{\bf Preuve.}  Pour simplifier les notations, \'ecrivons ${\cal O}_K$ pour ce
qui sera son localis\'e en ${\cal L}$, et $F$ son corps r\'esiduel (isomorphe
\`a $k({\cal L})$, donc dont la caract\'eristique est $l\neq p$). On a la
suite exacte de sch\'emas en groupes sur $F$ (comme faisceaux f.p.p.f.) :
$$0\to {\grdG}_{m,F} \to {\cal E}_F \to (\grdZ /N\grdZ )_F \to0.$$
Soit $a$ l'entier positif tel que ${\langle p^a .P \rangle}_F =\langle p^a .
\widetilde{P} \rangle =\langle \widetilde{P} \rangle \cap {\cal E}^0_F \,
(\subseteq  \langle p^k .\widetilde{P} \rangle ),$ o\`u on a not\'e
${\cal E}^0_F$ la composante neutre de la fibre sp\'eciale du mod\`ele de
N\'eron de $E$. On a $a\geq k$, et $p^a$ divise $N$. Il suffit de montrer que
${\cal E}_F /{\langle p^k .P\rangle}_F$ est un ouvert de ${\cal E}'_F$, car
alors le lemme est \'evident : on peut travailler dans ${\cal E}_F /{\langle
p^k .P\rangle}_F$, o\`u il est clair par le choix de $k$ que l'image de $P$
est d'ordre $p^k$ dans le groupe des composantes.\\  %
{\bf Premi\`ere \'etape.} Montrons qu'on a ${\cal E}_F /{\langle p^k .P
\rangle}_F ={({\cal E}/\langle p^k .P\rangle )}_F$. Le groupe fini $G:=(\grdZ
/p^{n-k} \grdZ )$ agit sur ${\cal E}_{/{\cal O}_K }$ par addition de $p^k .P$.
Comme ${\cal E}_{/{\cal O}_K }$ est quasi-projectif, l'orbite sous $G$ de
chaque point est contenue dans un affine. Donc on sait que tout ouvert affine
Spec($A$), stable par $G$, de ${\cal E}_{/{\cal O}_K}$ donne les ouverts
Spec($A^G {\otimes}_{{\cal O}_K} F$) et Spec($(A {\otimes}_{{\cal O}_K} F)^G
$) de ${({\cal E}/\langle p^k .P\rangle )}_F$ et ${\cal E}_F /{\langle p^k .P
\rangle}_F$ respectivement (voir \cite{serre}, III \S 12, \cite{mumford}, III
\S 12). On a donc juste \`a v\'erifier que le morphisme canonique entre les
deux anneaux ci-dessus est un isomorphisme. Il suffit pour cela de remarquer
que, $|G|$ \'etant inversible dans $A$, le projecteur $A\to A^G$ qui envoie
$x$ sur $\frac{1}{|G|} \sum_{G} g(x)$ commute au changement de base.\\
{\bf Deuxi\`eme \'etape.} Pour pouvoir d\'eduire de la propri\'et\'e
universelle des mod\`eles de N\'eron un prolongement \`a tout Spec(${\cal
O}_K$) du morphisme sur la fibre g\'en\'erique $({\cal E}/{\langle p^k .P
\rangle})_K \to {{\cal E}'}_K$, on doit v\'erifier que le premier sch\'ema est
lisse sur Spec(${\cal O}_K$). Consid\'erons la suite exacte de sch\'emas : %
$$
0 \to (\grdZ /p^{n-k}\grdZ)_{/{\cal O}_K} \to {\cal E}_{/{\cal O}_K}
\to ({\cal E}/{\langle p^k .P\rangle})_{/{\cal O}_K}\to 0 \ ;$$
ils sont tous ici localement de pr\'esentation finie, puisque localement de
type fini sur un anneau noeth\'erien. De plus, le premier sch\'ema de la suite
est \'etale sur ${\cal O}_K$ ; puisque $(\grdZ /p^{n-k} \grdZ )$ agit
librement sur ${\cal E}_{/{\cal O}_K}$, la seconde fl\`eche est finie \'etale
(voir \cite{Katz-Mazur}, th\'eor\`eme A7.1.1). \'Etant donn\'ees nos
hypoth\`eses, la lissit\'e sur ${\cal O}_K$ du dernier sch\'ema r\'esulte par
exemple de la proposition 17.7.7 de \cite{egaiv}.\\
{\bf Troisi\`eme \'etape.} On consid\`ere donc le morphisme prolong\'e $({\cal
E}/{\langle p^k .P\rangle})_{/{\cal O}_K} \to {{\cal E}'}_{/{\cal O}_K}$. Le
premier sch\'ema est s\'epar\'e sur la base. On peut donc
appliquer la proposition 3.2 de l'expos\'e IX de \cite{SGA7I} : et dire que ce
morphisme est une immersion ouverte. En se restreignant \`a la fibre
sp\'eciale, et en se servant de la premi\`ere \'etape, on peut bien voir comme
on le voulait plus haut ${\cal E}_F /{\langle p^k .P\rangle}_F$ comme un
ouvert de ${\cal E}'_F$. $\Box$
\section{L'alg\`ebre de Hecke pour ${\Gamma}_1$ en niveau et poids
quelconques.}
On a introduit dans le ``sch\'ema de la preuve'' l'alg\`ebre de Hecke sur
${\grdZ}$ pour ${\Gamma}_0$. On va maintenant se placer dans un cadre un peu
plus g\'en\'eral : ce qu'on notera dans cette section ${\grdT}_{\grdZ}$ ou
simplement ${\grdT}$ d\'esignera l'alg\`ebre de Hecke ``pour ${\Gamma}_1$, en
niveau et poids quelconques'' comme intitul\'e (pour plus de d\'etails, voir
la sous-section suivante). On %
notera ${\grdT}_{\grdQ}$ et ${\grdT}_{\grdC}$ les tensorisations de ${\grdT}$
avec ${\grdQ}$ et ${\grdC}$ respectivement. Le but de ce qui suit est
d'expliciter matriciellement l'action de ${\grdT}_{\grdQ}$ sur l'espace des
formes paraboliques de poids $\lambda$ pour ${\Gamma }_1 (N)$, o\`u $N$ est
un entier quelconque. Pour cela, on se sert des r\'esultats de
\cite{Atkin-Lehner} g\'en\'eralis\'es (voir par exemple \cite{Lang}, ou
\cite{Diamond} ; ces r\'ef\'erences seront constamment utilis\'ees). On en
d\'eduira la forme de ${\grdT}_{\grdQ}$ comme $\grdQ$-alg\`ebre abstrai\-te,
ce qui permettra ensuite, en se restreignant au cas de poids $2$ et de
caract\`ere trivial ({\it i.e.}, formes modulaires sur ${\Gamma}_0 (N)$), de
prouver la finitude du quotient d'enroulement sur $\grdQ$ en niveau quelconque
(th\'eor\`eme 3.9) - m\^eme si on n'a besoin pour le th\'eor\`eme 1.5 que des
niveaux puissance d'un nombre premier.
\subsection{Rappels.} %
On commence par fixer les notations : soit $f=\sum_{n\geq 1} a_n x^n$ une
s\'erie formelle  \`a coefficients dans un corps, en une variable ; soit
$\lambda$ et $N$ deux entiers, et $\varepsilon$ l'extension \`a $\grdZ$ d'un
caract\`ere de Dirichlet %
sur $(\grdZ /N \grdZ )^{\times}$ (on pose $\varepsilon (n)=0$ si $(n \wedge
N)\neq 1)$. On note $t_p$, $U_q$, $B_d$ ($p$, $q$, premiers, $d$ entier
quelconque), les op\'erateurs d\'efinis formellement (voir
\cite{Atkin-Lehner}) par :
$$\left\{ \begin{array}{l}
t_p (f) = \sum_{n\geq 1} (a_{np} +\varepsilon (p)p^{\lambda -1}\, a_{n/p} )
x^n \ ;\\
U_q (f)= \sum_{n\geq 1} a_{np} \, x^n \ ;\\
B_d (f)= \sum_{n\geq 1} a_{n} x^{nd} = \sum_{n\geq 1} a_{n/d} x^n \ ,
\end{array} \right.$$
o\`u on pose $a_{n/m} =0$ si $m$ ne divise pas $n$. On a les relations de
commutation suivantes :
$$\left\{ \begin{array}{l}
B_d \circ B_{d'} =B_{d'} \circ B_d ,{\mathrm{\ pour\ tous\ }}d,d'{\mathrm{\
dans\ }}\grdN \ ;\\
t_p \circ B_d =B_d \circ t_p ,{\mathrm{\ si\ }}p\ {\mathrm{et}}\ d\
{\mathrm{sont\ premiers\ entre\ eux}}\ ;\\
t_p \circ t_{p'}=t_{p'} \circ t_{p} {\mathrm{\ pour\ tous}}\  p\ et\ p'\
{\mathrm{premiers\ }};\\
t_p \circ U_q=U_q\circ t_p \ {\mathrm{si}}\ p\neq q\ ;\\
U_q\circ U_{q'}=U_{q'}\circ U_{q} \ ,\ {\mathrm{pour\ tous}}\ q\ {\mathrm{et}}
\ q'\ {\mathrm{premiers}}\ ;\\
U_q\circ B_d=B_d\circ U_q,\ {\mathrm{si}}\ d\ {\mathrm{et}}\ q'\ {\mathrm{sont
\ premiers\ entre\ eux.}}
\end{array} \right. $$
De plus,
$$U_q \circ B_{q^k} =B_{q^{k-1}} .$$
Consid\'erons maintenant l'espace des formes modulaires de poids $\lambda$
pour $\Gamma_1 (N)$ \`a coefficients dans $\grdZ$, $S_{\lambda} ({\Gamma}_1
(N))$, et sa tensorisation par $\grdQ$, $S_{\lambda} ({\Gamma}_1 (N))
{\otimes}_{\grdZ} \grdQ$, qu'on \'ecrira $S_{\lambda} ({\Gamma}_1 (N))_\grdQ$
(de m\^eme,  $S_{\lambda} ({\Gamma}_1 (N))_\grdC$). On peut d\'efinir
l'alg\`ebre de Hecke d'endomorphismes de ce $\grdQ$-espace, comme on le disait
pr\'ec\'edemment, qui est engendr\'ee par les op\'erateurs $T_p$ ($p$
premier), et les op\'erateurs diamants $\langle n\rangle$ (pour les $n$
premiers au niveau $N$). Un \'el\'ement de $S_{\lambda} ({\Gamma}_1 (N))_\grdQ
\subset S_{\lambda} ({\Gamma}_1 (N))_\grdC$ peut \^etre vu comme une fonction
$f(z)$ holomorphe sur le demi-plan de Poincar\'e, et on peut \'ecrire son
d\'eveloppement de Fourier en l'infini pour arriver \`a l'expres\-sion :%
$$f(x)=\sum_{n\geq 1} a_{n} \, x^n \ ,$$ %
o\`u $x=e^{2i\pi z}$, et les $a_n$ sont dans $\grdQ$. L'action des
op\'erateurs de Hecke $T_p$ sur  $S_{\lambda} ({\Gamma}_1 (N))_\grdQ$ est
alors pr\'ecis\'ement celle d\'ecrite plus haut par op\'erateurs $t_p$, si $p$
ne divise pas $N$, et
$U_q$ si $q$ divise $N$ ; on conserve d\'esormais cette notation ``mixte''
pour nous des $T_p \, ,U_q$, qui est celle d'Atkin-Lehner. On rappelle alors
le th\'eor\`eme principal de leur th\'eorie (voir \cite{Diamond},
\cite{Lang}) :
\begin{theo} %
{\bf (Atkin-Lehner)} L'espace $S_{\lambda} ({\Gamma}_1 (N))_\grdC$ se
d\'ecompose en une somme directe :%
$$S_{\lambda} ({\Gamma}_1 (N))_\grdC = \bigoplus_{\varepsilon} S_{\lambda}
(N, \varepsilon ),$$ %
o\`u $S_{\lambda} (N, \varepsilon )$ d\'esigne l'espace propre correspondant
au caract\`ere de Dirichlet $\varepsilon$, et la somme est prise sur tous les
tels $\varepsilon$ v\'erifiant $\varepsilon (-1)=(-1)^{\lambda }$.

  Chaque espace $S_{\lambda} (N, \varepsilon )_{\grdC}$ poss\`ede \`a son tour
une base ${\cal B}^{\varepsilon}$ compos\'ee de sous-bases ${\cal
B}_f^{\varepsilon}$ du type suivant :%
$${\cal B}_f^{\varepsilon} =\{ f(kz),\ k|(N/M)\} =\{ B_{k} (f),\ k|(N/M)\}
\ ,$$ %
o\`u $f$ est une newform en niveau $M$, de caract\`ere $\chi$ (sur $(\grdZ
/M\grdZ )^{\times}$) tel que l'extension de $\chi$ \`a $(\grdZ /N\grdZ
)^{\times}$ \'egale $\varepsilon$ ; si $C$ est le conducteur de $\varepsilon$
(ou de $\chi$), on a : $C|M|N$. Si $N=M$, alors ${\cal B}_f^{\varepsilon}$ ne
comporte qu'un \'el\'ement, et on l'appelle une {\em newclass} ; sinon, ${\cal
B}_f^{\varepsilon}$ est une {\em oldclass}. On note $E_f^{\varepsilon}$ le
sous $\grdC$-espace de $S_{\lambda} (N, \varepsilon )_{\grdC} \subset
S_{\lambda} ({\Gamma}_1 (N))_\grdC$ engendr\'e par ${\cal B}_f^{\varepsilon}$.

  Pour $p$ premier \`a $M$, (respectivement, $q$ premier divisant $M$), $f$
est un vecteur propre de $T_p$ (respectivement $U_q$), et la valeur propre
associ\'ee est le $p$-i\`eme coefficient $a_p$ de $f$, (qui, puisque newform,
est suppos\'ee normalis\'ee par $a_1 =1$).

  Dans le cas o\`u $\varepsilon$ est trivial (formes modulaires pour
${\Gamma}_0 (N)$), si $q^2$ divise $M$, on a $U_q (f)=0$, tandis que si $q$ -
mais pas son carr\'e - divise $M$ (ce qu'on note $q||M$), $U_q (f)=\pm f$.
\end{theo}%
On va donc maintenant \'ecrire matriciellement l'action de l'alg\`ebre de
Hecke ${\grdT }_\grdC$ sur $S_{\lambda} ({\Gamma}_1 (N))_\grdC$, en
\'ecrivant cet espace comme somme de facteurs $S_{\lambda} (N, \varepsilon
)_{\grdC}$, qu'on d\'ecompose \`a leur tour en somme des $E_f^{\varepsilon}$
qui correspondent aux classes du th\'eor\`eme. (On en d\'eduira \`a chaque
fois l'action de ${\grdT}_{\grdQ}$ sur $S_{\lambda} ({\Gamma}_1 (N))_\grdQ$,
en le d\'ecomposant en sous-$\grdQ$-espaces ${\cal E}_f^{\varepsilon}$ dont la
tensorisation avec $\grdC$ donne la somme des conjugu\'es sous Galois des
$E_f^{\varepsilon}$.) On commencera d'abord (3.3) par \'etudier le cas de
``co-niveau'' puissance de premier, c'est-\`a-dire le cas o\`u, avec les
notations du th\'eor\`eme, $(N/M)=p^k$ pour $p$ premier - (et c'est encore une
fois le seul cas dont on ait besoin dans le reste du papier). Puis on traitera
le cas g\'en\'eral en 3.4.
\subsection{Cas de co-niveau qu'un seul premier $p$ divise.}
Soit donc $f\in S_{\lambda} (N, \varepsilon )_{\grdC}$, qui est une newform
en niveau $M$, avec $N/M=~p^k$. On a alors ${\cal B}_f =\{ (f(z),\,  f(pz),\,
f(p^2 z),\, \dots ,f(p^{k} z)\} = \{ f, B_p (f),\, \dots ,B_{p^k} (f)\}$.
Examinons l'action des diff\'erents op\'erateurs de Hecke sur $E_f$ :
d'apr\`es les r\'esul\-tats pr\'ec\'edents, les op\'erateurs $T_l$, $l$ ne
divisant pas $N$, et $U_q$, $q$ divisant $M$ et diff\'erent de $p$, de m\^eme
que les op\'erateurs diamants, agissent diagonalement sur $E_f$, puisqu'ils
ont $f$ comme vecteur propre et qu'ils commutent avec les $B_{p^j} ,\ 0\leq j
\leq k$. Puisque ${\grdT}_{\grdZ}$ op\`ere sur $S_{\lambda} ({\Gamma}_1 (N)
)_\grdZ ,$ qui est un $\grdZ$-module libre de rang fini, les valeurs propres
des $T_l$ sont des entiers alg\'ebriques. Et il existe un corps de nombre
$K_f$ contenant tous les $a_p$ (car l'alg\`ebre de Hecke est un $\grdZ$-module
de type fini).\\
{\bf Remarque.} {\it Le corps de nombres engendr\'e par les valeurs propres
associ\'ees \`a $f$ de ${\grdT}_{\grdQ}$ n'est pas plus grand que $K_f$ : en
effet, m\^eme s'il contient, en plus des $a_n$, les valeurs propres des
op\'erateurs diamants, la relation $p^{\lambda -1}\langle p\rangle =T_p^2 -
T_{p^2}$ (pour tout $p$ premier ne divisant pas le niveau) montre que ces
valeurs propres appartiennent \`a $K_f$.}

  En caract\`ere trivial, les $U_q$ sont triviaux, puisque leurs valeurs
propres en $f$ sont soit $0$, soit $1$, soit $-1$ ; dans ce cas encore, les
coefficients $a_n$ de $f$ sont de plus totalement r\'eels : en effet, les
op\'erateurs de Hecke sont auto-adjoints pour le produit scalaire de
Petersson :%
$$\langle f,g\rangle =\int \int_{z =x+iy\, \in D} f(z)\, \overline{g(z)} \,
y^{\lambda -2} \, dx\, dy\ ,$$%
qui fait de $S_2 (N)_\grdC$ un espace de Hilbert (on a d\'esign\'e par $D$
dans l'int\'egrale un domaine fondamental du demi-plan de Poincar\'e pour
$\Gamma_0 (N)$). (En fait, on d\'efinit un produit scalaire de Petersson pour
tout $S_{\lambda} ({\Gamma}_1 (N))_\grdC$, pour lequel les $T_l$ sont normaux
(sinon auto-adjoints), et pour lequel la d\'ecomposition en
$E_f^{\varepsilon}$ est orthogonale - voir la preuve de 3.7.)
\setlength{\unitlength}{0.7cm}

  Revenant au cas g\'en\'eral, on consid\`ere maintenant l'action de $U_p$.
Pour cela, on distingue deux cas : selon que $p$ divise ou non $M$.
\subsubsection{Si $p$ divise $M$.} %
L'op\'erateur $U_p$ de l'alg\`ebre de Hecke en niveau $N$ est le m\^eme que
celui de niveau $M$, et donc avec le lemme plus haut,%
$$\left\{ \begin{array}{l}
U_p (f)=a_p \, f,\ {\mathrm et}\\
U_p (B_{p^j} (f))=B_{p^{j-1}} (f)\ {\mathrm si}\ j\geq 1\ .
\end{array} \right.$$ %
Dans la base ${\cal B}_f$, $U_p$ est donc sous la forme $(M_1 )$ :
$$U_p = \left(
\begin{picture}(7,4)(0.6,4)
\put(1,7){$a_p$}
\put(2,7){1}
\put(1,6){0}
\put(2,6){0}
\put(3,6){1}
\multiput(3,7)(1,0){5}{.}
\put(7,7){0}
\multiput(1,6)(0,-1){5}{.}
\put(1,1){0}
\multiput(2,1)(1,0){5}{.}
\put(7,1){0}
\multiput(7,6)(0,-1){4}{.}
\put(7,2){1}
\multiput(3,5)(1,-1){4}{.}
\multiput(4,5)(1,-1){3}{.}
\end{picture}
\right) .$$%
Si $a_p =0$, la matrice pr\'ec\'edente n'a qu'un bloc de Jordan, et la
restriction $R_f$ de l'alg\`ebre de Hecke ${\grdT}_{\grdQ}$ au $\grdQ$-espace
${\cal E}_f^{\varepsilon}$ qu'on a d\'efini plus haut (correspondant sur
$\grdC$ \`a la somme des conjugu\'es par Galois de $E_f^{\varepsilon}$), est
de forme :%
$$R_f =K_f [U_p ]\simeq K_f [X]/(X^{k+1} ).$$ %
Sinon, on a pour $U_p$ deux blocs de Jordan, un petit correspondant
\`a la valeur propre $a_p$, et un gros, correspondant \`a $0$ ; en tant que
$\grdQ$-alg\`ebres, on a donc l'isomorphisme
$$R_f \simeq K_f \times K_f [X]/(X^k ).$$ %
{\bf Remarque.} {\it Le fait que le corps du deuxi\`eme facteur de $R_f$ soit
bien tout $K_f$, malgr\'e l'absence de la valeur propre $a_p$, peut se voir
par un argument de dimension des espaces cotangents ; ou bien avec les
r\'esultats d'Atkin-Lehner, qui disent moralement qu'une newform est
caract\'eris\'ee par ``tous ses coefficients de Fourier moins un nombre
fini''.} %
\paragraph{Application au cas de caract\`ere trivial.} %
Dans ce cas, si $p^2$ divise $M$, $a_p =0$, et $R_f \simeq K_f [X]/(X^{k+1})$
; la base ${\cal B}_f^{\varepsilon}$ diagonalise $U_p$. Si en revanche
$p$ - mais pas son carr\'e - divise $M$, on a $a_p =\pm 1$, donc $R_f \simeq
K_f \times K_f [X]/(X^k )$, et on voit facilement que la nouvelle base de
$E_f^{\varepsilon}$ :%
$${\cal B}'_f=\{ f,\, B_p (f)-a_p f,\, B_{p^2} (f)-f,\, \ldots ,\,
B_{p^{2j}} (f)-f,\,  B_{p^{2j+1}} (f)-a_p  f,\ldots \}$$ %
est de Jordan pour $U_p$ et donc triangulise toute la restriction de
l'alg\`ebre de Hecke \`a $E_f^{\varepsilon}$.
\subsubsection{Si $p$ ne divise pas $M$.}
Le fait nouveau par rapport au cas pr\'ec\'edent est que l'alg\`ebre de
Hecke en niveau $N$, ${\grdT}_{N,\grdC}$, qu'on fait agir sur
$E_f^{\varepsilon}$, n'est plus la restriction \`a cet espace de l'alg\`ebre
de Hecke en niveau $M$ : le ``bon'' $p$-i\`eme op\'erateur de Hecke, celui
pour lequel $f$ est un vecteur propre, est $T_p$, mais puisqu'on s'int\'eresse
\`a $\grdT$ en niveau $N$, on doit consid\'erer \`a la place l'op\'erateur
$U_p$.

  Puisque $f$ est propre pour $T_p$, on a :%
$$T_p (f)=\sum_{n\geq 1} (a_{np} +\varepsilon (p)\, p^{\lambda -1}\, a_{n/p} )
\, x^n \, =a_p \, f\ ,$$ %
et les coefficients de $f$ v\'erifient donc pour tout $n$ les relations :%
$$a_p \, a_n = a_{np} +\varepsilon (p)\, p^{\lambda -1}\, a_{n/p} .$$ %
On calcule alors :%
\begin{eqnarray*}
U_p (f) & = & \sum_{n\geq 1} a_{np} \, x^n \\
 & = & \sum_{n\geq 1} (a_p a_n \, -\varepsilon (p)\, p^{\lambda -1} \, a_{n/p}
)\, x^n \\
 & = & a_p f -\varepsilon (p)\, p^{\lambda -1} \, B_{p} (f).
\end{eqnarray*}
La matrice de $U_p$ dans ${\cal B}_f$ est donc $(M_2 )$ :
$$U_p = \left(
\begin{picture}(9,4)(0.6,4)
\put(2,7){$a_p$}
\put(.6,6){$(-\varepsilon (p)\, p^{\lambda -1} )$}
\put(2,5){0}
\put(4,7){1}
\put(4,6){0}
\put(5,6){1}
\put(5,5){0}
\put(6,5){1}
\multiput(5,7)(1,0){5}{.}
\put(9,7){0}
\multiput(2,4)(0,-1){3}{.}
\put(2,1){0}
\multiput(4,1)(1,0){5}{.}
\put(9,1){0}
\multiput(9,6)(0,-1){4}{.}
\put(9,2){1}
\multiput(6,4)(0.5,-0.5){4}{.}
\multiput(7,4)(0.5,-0.5){3}{.}
\end{picture}
\right) .$$
Le polyn\^ome caract\'eristique de la restriction de $U_p$ \`a
$E_f^{\varepsilon}$ est %
$${P}_{U_p} (X)=(X^2 -a_p X+\varepsilon (p)\, p^{\lambda -1}).X^{k-1} \ ,$$
et les racines $\alpha_p$ et $\overline{\alpha}_p$ du premier facteur sont
simples si ${a_p }^2 \neq 4\varepsilon (p)\, p^{\lambda -1}$, doubles (valant
$\pm \sqrt{\varepsilon (p)} p^{\frac{\lambda -1}{2}}$) sinon ; dans
ce deuxi\`eme cas on voit sur le rang de $U_p$ qu'il a deux blocs de Jordan,
un pour chaque valeur propre ($0$ et $\pm \sqrt{\varepsilon (p)}
p^{\frac{\lambda -1}{2}}$).\\  %
{\bf Remarque.} {\it En caract\`ere trivial et poids $2$, ce dernier cas, (qui
correspond \`a une valeur maximale de $a_p$ selon la borne de Weil, qui dit
que $|a_p |\leq 2\sqrt{p}$), est en fait exclu, par un th\'eor\`eme non encore
publi\'e de Coleman et Edixhoven (\cite{Coleman-Edixhoven}) ; on l'explicite
quand m\^eme - en attendant.} \\  %
\underline{Si ${a_p}^2 \neq 4\varepsilon (p)\, p^{\lambda -1}$,} dans une
base convenable la restriction de $U_p$ \`a $E_f^{\varepsilon}$ s'\'ecrit
$(M_3 )$ : %
$$U_p = \left(
\begin{picture}(8,4.5)(-0.4,3.5)

\put(-0.5,5.7){\dashbox{.3}(2,1.8){ }}
\put(0.7,6){$\overline{\alpha}_p$}
\put(1,7){0}
\put(0,6){0}
\put(0,7){$\alpha_p$}

\put(0.5,3){0}

\put(4.5,6.5){0}

\put(1.6,-0.4){\dashbox{.4}(6,6){ }}
\put(2,5){0}
\put(3,5){1}
\put(3,4){0}
\put(4,4){1}
\multiput(4,5)(1,0){3}{.}
\put(7,5){0}
\multiput(2,4)(0,-1){4}{.}
\put(2,0){0}
\multiput(3,0)(1,0){4}{.}
\put(7,0){0}
\multiput(7,4)(0,-1){3}{.}
\put(7,1){1}
\multiput(4,3)(1,-1){3}{.}
\multiput(5,3)(1,-1){2}{.}
\end{picture}
\right) .$$
Dans ce cas-l\`a, si on pose $K'_f =K_f [X]/(X^2 -a_p X+\varepsilon (p)\,
p^{\lambda -1} )$, (qui n'est pas n\'ec\'essairement un corps), la restriction
de ${\grdT}_{\grdQ}$ \`a ${\cal E}_f^{\varepsilon}$ est %
$$R_f \simeq K'_f \times K_f [X]/(X^{k-1} ).$$  %
On note (cela servira dans la suite) qu'un vecteur propre associ\'e \`a la
valeur propre $0$ est $(B_{p^2} (f)-\frac{a_p}{\varepsilon (p)\,p^{\lambda
-1}} B_p (f)+\frac{1}{\varepsilon (p)\, p^{\lambda -1}} f)$. \\
\underline{Si ${a_p}^2 =4\varepsilon (p)\, p^{\lambda -1},$} on a de m\^eme
dans une base convenable l'\'ecriture matricielle $(M_4 )$ : %
$${U_p}_{|E_f^{\varepsilon}} = \left(
\begin{picture}(8,4.5)(-0.4,3.5)
\put(-0.5,5.7){\dashbox{.3}(2,1.8){ }}
\put(0.7,6){$\frac{a_p}{2}$}
\put(1,7){1}
\put(0,6){0}
\put(0,7){$\frac{a_p}{2}$}
\put(0.5,3){0}
\put(4.5,6.5){0}
\put(1.6,-0.4){\dashbox{.4}(6,6){ }}
\put(2,5){0}
\put(3,5){1}
\put(3,4){0}
\put(4,4){1}
\multiput(4,5)(1,0){3}{.}
\put(7,5){0}
\multiput(2,4)(0,-1){4}{.}
\put(2,0){0}
\multiput(3,0)(1,0){4}{.}
\put(7,0){0}
\multiput(7,4)(0,-1){3}{.}
\put(7,1){1}
\multiput(4,3)(1,-1){3}{.}
\multiput(5,3)(1,-1){2}{.}
\end{picture}
\right) .$$ %
Alors avec les notations pr\'ec\'edentes,
$$R_f \simeq K_f [X]/(X^2 )\times K_f [Y]/(Y^{k-1} ).$$
Une base de Ker$(U_p -\alpha_p$ Id) est $\{ a_p f-2\varepsilon (p)p^{\lambda
-1} B_p (f)\}$, et bien s\^ur le m\^eme vecteur $(B_{p^2} (f)-
\frac{a_p}{\varepsilon (p)\, p^{\lambda -1}} B_p (f)+\frac{1}{\varepsilon (p)
\, p^{\lambda -1}} f)$ que pr\'ec\'edemment engendre Ker($U_p$).
\subsection{Cas de co-niveau quelconque.} %
Soit maintenant $f\in  S_{\lambda} (N, \varepsilon )_{\grdC}$, qui est une
newform en niveau $M$, avec $N/M=n$ quelconque cette fois. On a alors
${\cal B}_f =\{ B_d (f), d|n\}$. On examine \`a nouveau l'action des
diff\'erents op\'erateurs de Hecke : les op\'erateurs $T_l$, $l$ ne divisant
pas $N$, et $U_q$, $q$ ne divisant pas $n$, agissent toujours diagonalement
sur $E_f^{\varepsilon}$, comme en co-niveau premier. Soit $q$ un premier
divisant $n$. Notons $m$ la valuation en $q$ de $n$ : $q^m ||n$. Alors, on
construit pour $U_q$ une base qui lui est appropri\'ee, en ordonnant
partiellement la base ${\cal B}_f$ ainsi : pour $d$ parcourant les diviseurs
de $n/{q^m}$, on note ${\cal B}_f^{q,d} =\{ B_{d} (f),\, B_{dq} (f),\,
B_{dq^2} (f),\, \ldots ,\, B_{dq^m} (f)\}$, et on \'ecrit %
$${\cal B}_f^q =\{ f, B_q (f), B_{q^2} (f), \ldots ,B_{q^m} (f), B_{d} (f),
B_{dq} (f), B_{dq^2} (f), \ldots , B_{dq^m} (f), \ldots \}$$
comme la r\'eunion sur $d$ de ces ${\cal B}_f^{q,d}$. (Le nombre des
diff\'erentes bases ${\cal B}_f^{q,d}$ est ${\sigma}_0 (n/{q^m})$ (nombre de
diviseurs de $n/{q^m}$).) On d\'eduit des relations de commutation entre les
op\'erateurs consid\'er\'es et de ce qui pr\'ec\`ede que, dans la base
${\cal B}_f^{q}$, on a :
$${U_q}|_{E_f} =\left( \begin{array}{cccc}
B^q_1 & 0 & \cdots & 0 \\
0 & B^q_{d_2} &  \cdots & 0 \\
\vdots &  & \ddots &  \\
0 & \ldots & \ldots & B^q_{d_{{\sigma}_0 (n/{q^m})}}
\end{array} \right) \ ,$$
o\`u chaque $B^q_d,\ d|(n/{q^m})$ est une matrice repr\'esentant
la restriction de $U_q$ \`a l'espace engendr\'e par  ${\cal B}_f^{q,d}$ - que
$U_q$ laisse stable. La forme de $B^q_d$ ne d\'epend en fait pas de $d$, mais
uniquement de $q$. En effet :\\
\underline{Si $q$ divise $M$}, alors chaque $B^q_d$ est une bloc de forme
$(M_1 )$, et on peut trouver une nouvelle base de l'espace engendr\'e par
${\cal B}_f^{q,d}$, dans laquelle on a $B^q_d$ \'equivalente si $a_p \neq 0$
\`a %
$$\left(
\begin{picture}(7,4)(0.6,4)

\put(0.5,6.7){\dashbox{.3}(0.8,0.8){$a_q$}}

\put(1,3.5){0}

\put(4.5,7){0}

\put(1.6,0.6){\dashbox{.4}(6,6){ }}

\put(2,6){0}
\put(3,6){1}
\put(3,5){0}
\put(4,5){1}
\multiput(4,6)(1,0){3}{.}
\put(7,6){0}
\multiput(2,5)(0,-1){4}{.}
\put(2,1){0}
\multiput(3,1)(1,0){4}{.}
\put(7,1){0}
\multiput(7,5)(0,-1){3}{.}
\put(7,2){1}
\multiput(4,4)(1,-1){3}{.}
\multiput(5,4)(1,-1){2}{.}
\end{picture}
\right) ,$$
et sinon \`a un bloc de Jordan nilpotent.\\  %
\underline{Si $q$ ne divise pas $M$}, on conclut toujours comme dans la
sous-section pr\'ec\`edente que $B^q_d =M_2 $, \'equivalente si ${a_q}^2
\neq 4\varepsilon (q) q^{\lambda -1}$ \`a $(M_3 )$, et si ${a_q}^2  =4
\varepsilon (q)q^{\lambda -1}$ \`a $(M_4 )$. Maintenant, si on note $E_f^q$ le
sous-espace de $S_{\lambda} (N, \varepsilon )_{\grdC}$ engendr\'e par
${\cal B}_f^{q,1} =\{ f,B_{q} (f),\, B_{q^2} (f),\, \ldots ,\, B_{q^m} (f)\}$,
alors $U_q$ agit comme pr\'ecis\'e ci-dessus sur $E_f^q$, et on peut
consid\'erer $E_f^{\varepsilon}$ comme le produit tensoriel sur $\grdC$ : %
$$E_f^{q_1} \otimes  E_f^{q_2} \otimes \cdots \otimes E_f^{q_m} ,$$ %
o\`u les $q_i$ sont les nombres premiers divisant $n$. Et $S_{\lambda} (N,
\varepsilon )_{\grdC}$, \`a son tour, est la somme directe des tels
$E_f^{\varepsilon}$.

  On r\'esume tout \c{c}a :%
\begin{theo}
L'entier naturel $N$ est quelconque.

   L'alg\`ebre de Hecke ${\grdT}_{\grdQ}$, vue comme sous-alg\`ebre des
endomorphismes de $\grdQ$-espace vectoriel de $S_{\lambda} ({\Gamma}_1
(N))_\grdQ$, s'\'ecrit comme un produit d'anneaux :%
$${\grdT}_{\grdQ} =R_{f_1} \times R_{f_2} \times \cdots \times
R_{f_k} ,$$ %
o\`u chacun de ces $R_{f_i}$ correspond \`a l'orbite sous Galois d'une newform
$f_i$ en niveau $M_i$ divisant $N$, et caract\`ere $\varepsilon$.

  Plus pr\'ecis\'ement, le facteur $R_{f_i}$ est la restriction de
${\grdT}_{\grdQ}$ au sous-espace ${\cal E}_{f_i}$ de $S_{\lambda} ({\Gamma}_1
(N))_{\grdQ}$, qui est engendr\'e sur $\grdQ$ par les orbites sous Galois des
$B_d (f_i )$, $d$ parcourant les diviseurs de $N/M_i$.

  Si $K_{f_i}$ d\'esigne le corps de nombres engendr\'e par les coefficients
$a_n$ du d\'eveloppement de Fourier de $f_i$ \`a l'infini, l'id\'eal $R_{f_i}$
est de forme un produit tensoriel d'anneaux des quatre types suivants :
$$\begin{array}{l}
\bullet \ K_{f_i} [X]/(X^{n} )\ ;\\
\bullet \ K_{f_i} \times K_{f_i} [X]/(X^n )\ ;\\
\bullet \ K_{f_i} [X]/(X^2 -a_q X+\varepsilon (q)q^{\lambda -1} )\times
K_{f_i} [Y]/(Y^n )\ ;\\
ou\ encore\ :\\
\bullet \ K_{f_i} [X]/(X^2 )\times K_{f_i} [Y]/(Y^n ).
\end{array}$$
\end{theo}
Notons un corollaire partiel de ce th\'eor\`eme qui nous servira dans la
suite :
\begin{cor}
Soit $A$ une alg\`ebre de dimension finie sur un corps $K$, et qui a un unique
id\'eal minimal non trivial ; alors $A$ est de Gorenstein, i.e. son dual en
tant que $K$-espace vectoriel, $A^{\vee}$, est un $A$-module libre de
rang 1. En particulier, l'alg\`ebre de Hecke sur $\grdQ$ (pour ${\Gamma}_1$,
en niveau et poids quelconques) ${\grdT}_{\grdQ}$, est de Gorenstein.
\end{cor}
{\bf Preuve.}\ Soit $A$ une alg\`ebre comme dans l'\'enonc\'e. Son id\'eal
minimal est par d\'efinition principal, engendr\'e par n'importe lequel de ses
\'el\'ements non nuls. Soit $g$ l'un de ceux-l\`a. Il existe un \'el\'ement
$L$ de $A^{\vee}$ qui ne s'annule pas en $g$. On en d\'eduit donc que
le noyau de l'application lin\'eaire :
$$\left\{ \begin{array}{ccl}
A & \to & A^{\vee} \\
a & \mapsto & (x\mapsto L(a.x))
\end{array} \right.$$
est un id\'eal trivial de $A$, puisqu'il ne contient pas l'id\'eal minimal.
Donc cette application est injective, et en fait bijective puisque les deux
espaces $A$ et $A^{\vee}$ ont m\^eme dimension.

  Pour ce qui est de ${\grdT}_{\grdQ}$, le th\'eor\`eme pr\'ec\'edent dit
qu'elle se d\'ecompose en produits d'anneaux de la forme $K[X_1 ,X_2 ,\dots ,
X_n ]/(X_1^{m_1} X_2^{m_2} \dots X_n^{m_n} )$, o\`u $K$ est un corps de
nombres (en effet, chaque $A_{f_i}$ se d\'ecompose encore ainsi). Puisque
$K[X_1 ,X_2 ,\dots ,X_n ]/(X_1^{m_1} X_2^{m_2} \dots X_n^{m_n} )$ poss\`ede un
seul id\'eal minimal, %
et puisque la propri\'et\'e d'\^etre de Gorenstein est clairement stable par
somme directe finie, on peut appliquer la premi\`ere partie de la proposition.
$\Box$

  Montrons aussi un r\'eultat dont on se sert dans la sous-section suivante
(et qui sera red\'emontr\'e en partie en 4.7 gr\^ace \`a 3.6) :
\begin{prop}
Pour tous niveau $N$ et poids $\lambda$, $S_{\lambda} ({\Gamma}_1
(N))_{\grdQ}$ est un ${\grdT}_{\grdQ}$-module libre de rang 1.
\end{prop}
{\bf Preuve}.\ On a un accouplement parfait $[\ ,\ ]$, $\grdC$-bilin\'eaire :
$$\left\{ \begin{array}{rcl}
S_{\lambda} ({\Gamma}_1 (N))_{\grdC} \times {\grdT}_{\grdC} & \to & \grdC \\
(f=\sum_{n\geq 1} a_n q^n ,T) & \mapsto & a_1 (Tf) \ ,
\end{array} \right.$$
qui est d\'efini sur ${\grdQ}$, pour lequel tout op\'erateur de l'alg\`ebre de
Hecke est auto-adjoint. Maintenant, notons $\langle \  ,\,
\rangle_{\mathrm P}$ le produit scalaire de Petersson. Appelons %
$\sigma$ l'involution sur $S_{\lambda} ({\Gamma}_1 (N))_{\grdC}$ d\'efinie par
l'op\'eration de conjugaison complexe des coefficients de Fourier en l'infini
d'une forme parabolique. Rappelons aussi l'op\'erateur (classique) $W_N$,
dont l'op\'eration sur $S_{\lambda} ({\Gamma}_1 (N))_{\grdC}$ est d\'efinie
par $W_N (f(z))=N^{-\lambda /2} z^{-\lambda} f(-1/Nz)$ ; il v\'erifie
$W_N^2 =(-1)^{\lambda}$, et l'adjoint pour le produit scalaire de Petersson
de tout op\'erateur de Hecke $T_n$ est $W_N T_n W_N^{-1}$ (de m\^eme,
l'adjoint de l'op\'erateur diamant $\langle n \rangle$ pour $n$ premier \`a
$N$ est $W_N \langle n\rangle W_N^{-1}$) (voir par exemple \cite{Diamond},
I.4). Consid\'erons alors l'alt\'eration suivante : $(f,g)\to \{ f,g\} =
\langle f,\sigma \circ W_N (g)\rangle_{\mathrm P}$ du produit scalaire de
Petersson, c'est-\`a-dire :
$$\left\{ \begin{array}{rcl}
S_{\lambda} ({\Gamma}_1 (N))_{\grdC} \times S_{\lambda} ({\Gamma}_1
(N))_{\grdC} & \to & \grdC \\
(f,g) & \mapsto & \int_{\Delta} (\sum_{n\geq 1} a_n q^n )\, W_N
(\sum_{n\geq 1} b_n {\overline q}^n )\, y^{\lambda} \frac{dx\, dy}{y^2} \ ,
\end{array} \right.$$
o\`u $(\sum_{n\geq 1} a_n q^n )$ et $(\sum_{n\geq 1} b_n q^n )$ sont les
d\'eveloppements de Fourier de $f$ et $g$ respectivement en l'infini, et
$\Delta$ est un domaine fondamental du quotient du demi-plan de Poincar\'e par
${\Gamma}_1 (N)$. Par construction, ce produit est $\grdC$-bilin\'eaire, et
pour lui \'egalement l'alg\`ebre de Hecke est auto-adjointe. Les deux
accouplements $[\ ,\ ]$ et $\{ \ ,\ \}$ pr\'esentent donc respectivement
${\grdT}_{\grdC}$ et $S_{\lambda} ({\Gamma}_1 (N))_{\grdC}$ comme les duaux
comme $\grdC$-espaces de $S_{\lambda} ({\Gamma}_1 (N))_{\grdC}$, ce qui donne
un isomorphis\-me de $\grdC$-espaces vectoriels entre eux. Mais
l'auto-adjonction de ${\grdT}_{\grdC}$ pour les deux montre que cette
application est aussi un morphisme de ${\grdT}_{\grdC}$-module, donc que
$S_{\lambda} ({\Gamma}_1 (N))_{\grdC}$ est un ${\grdT}_{\grdC}$-module libre
de rang 1. Enfin, que ceci soit vrai sur ${\grdQ}$ est une cons\'equence du
r\'esultat de g\'eom\'etrie alg\'ebrique \'el\'ementaire suivant : si $k$ est
un corps infini dont $K$ est une extension, si $R$ est une $k$-alg\`ebre et
$M$ un $R$-module qui est un $k$-espace vectoriel de dimension finie, et tel
que $M\otimes_{k} K$ soit un $R\otimes_{k} K$-module libre de rang $1$, alors
$M$ est soi-m\^eme un $R$-module libre de rang $1$. Ce r\'esultat se montre
comme suit. Choisisons des \'el\'ements $r_i$ de $R$, $1\leq i\leq {\mathrm
{dim}}_k M$, tel que les $r_i .y$ forment une $K$-base de $M\otimes_{k} K$
pour un $y$ de $M\otimes_{k} K$. La fonction $f$ qui \`a tout \'el\'ement $x$
de $M$ associe le d\'eterminant de $(r_i .x)$ est polyn\^omiale en les
coordonn\'ees de $x$ dans une $k$-base de $M$, et n'est pas nulle sur
$M\otimes_{k} K$. Ce qui veut dire que le polyn\^ome correspondant est non
nul, et comme $k$ est infini, on a bijection entre polyn\^omes et fonctions
polyn\^omiales sur $M$ : donc $f$ ne peut \^etre uniform\'ement nulle sur $M$,
et $M$ est un $R$-module de rang 1. Qu'enfin $M$ soit $R$-libre provient de sa
$R\otimes_{k} K$-fid\'elit\'e apr\`es extension des scalaires \`a $K$. $\Box$
\subsection{Finitude du quotient d'enroulement.}
On la d\'emontre en niveau quelconque, avec les r\'esultats qui pr\'ec\`edent
appliqu\'es au cas de poids $\lambda$ \'egal \`a $2$, le caract\`ere
$\varepsilon$ \'etant trivial. Soit $N$ un entier ; on consid\`ere la courbe
modulaire $X_0 (N)_{\grdQ}$ et sa jacobienne $J_0 (N)_{\grdQ}$ sur $\grdQ$. On
rappelle les notations de l'introduction, dans ce cadre plus g\'en\'eral : $e$
d\'esigne l'\'el\'ement d'enroulement, c'est-\`a-dire l'\'el\'ement de $H_1
(X_0 (N)\, ;\, \grdQ )$ d\'efini par l'int\'egration entre z\'ero et l'infini
sur notre courbe modulaire ; ${\cal A}_e$ d\'esigne l'id\'eal d'enroulement,
c'est-\`a-dire l'id\'{e}al annulateur dans ${\grdT}_{\grdZ}$ de~$e$. Notons
aussi ${\cal A}_{e,\grdQ} ={\cal A}_e \otimes_{\grdZ} \grdQ$. Le quotient
d'enroulement  $J_0^e$ est la vari\'et\'e ab\'{e}lienne quotient $J_0 (N)/
{\cal A}_e J_0 (N)$. %
\begin{theo}
$J_0^e (\grdQ )$ est fini.
\end{theo}
{\bf Preuve.}\ \ Notons d'abord $S_2 (N)_{\grdQ}$ le $\grdQ$-espace des formes
paraboliques de poids deux pour $\Gamma _0 (N)$ \`a coefficients dans $\grdQ$,
et ${\mathrm {Cot}}_0 (J_0 (N)_{\grdQ} )$ l'espace cotangent \`a $J_0 (N
)_{\grdQ}$ en z\'ero ; on a les identifications suivantes :%
$$S_2 (N)_{\grdQ} \simeq H^0 (X_0 (N)_{\grdQ} , \Omega^1 )\simeq H^0
(J_0 (N)_{\grdQ} , \Omega^1 )\simeq {\mathrm {Cot}}_0 (J_0 (N)_{\grdQ} ),$$
et  $S_2 (N)_{\grdQ}$ est un $\grdT {\otimes }_{\grdZ} \grdQ$-module
libre de rang 1 (3.7, ou la proposition 4.7 de la section suivante).

  L'alg\`ebre de Hecke \`a coefficients dans $\grdQ$ s'\'ecrit donc comme on
l'a explicit\'e plus haut, {\em i.e.} comme un produit d'anneaux (ou comme une
somme directe d'id\'eaux) $R_f$, correspondants aux restrictions de
${\grdT}_{\grdQ}$ aux orbites sous l'action du groupe de Galois des newforms
$f$ de niveau $M$ divisant $N$.

   D'apr\`es le th\'eor\`eme d'Atkin-Lehner (th\'eor\`eme 3.2), si on note
$S_2 (M)^{\mathrm {new}} _{\grdC}$ le sous-$\grdC$-espace de  $S_2
(M)_{\grdC}$ engendr\'e par les {\em newforms} de niveau $M$, on a en faisant
la somme directe des op\'erateurs $B_d$ convenables un isomorphisme de
$\grdC$-espaces vectoriels :
$${\oplus}_{M|N} {\oplus}_{d|(N/M)} S_2 (M)^{\mathrm {new}} _{\grdC}
\stackrel{\sim}{\to} S_2 (N)_{\grdC} .$$
Gr\^ace \`a l'interpr\'etation en termes d'espaces cotangents des espaces de
formes paraboliques, on d\'eduit des op\'erateurs $B_d$ des isog\'enies
$J_0 (M)_{\grdQ} \stackrel{B^*_d}{\rightarrow} J_0 (N)_{\grdQ}$ (les
op\'erateurs $B_d$ sont d\'efinis sur ${\grdQ}$), et pour tout $M$ on
d\'efinit $J_0 (M)_{\mathrm {new}} :=  J_0 (M)/(\sum_{d|M}
{\mathrm {Im}}({B^*_d}))$. On d\'eduit alors de l'isomorphisme pr\'ec\'edent
(entre espaces de formes paraboliques) l'isog\'enie :
$$J_0 (N)_{\grdQ} \to {\oplus}_{M|N} {\oplus}_{d|(N/M)}
J_0 (M)_{\grdQ ,{\mathrm {new}}}  .$$
On a, encore et toujours selon la d\'ecomposition d'Atkin-Lehner, des
isog\'enies $ J_0 (M)_{\grdQ ,{\mathrm {new}}} \to {\oplus}_{G_{\grdQ} f}
J_f$, o\`u la somme est prise sur les orbites sous Galois $G_{\grdQ} f$ des
newforms $f$ en niveau $M$, et $J_f$ est la vari\'et\'e ab\'elienne
``d\'ecoup\'ee'' dans  $J_0 (M)_{\grdQ}$ avec la forme $f$ : avec les
notations de 3.5, l'annulateur ${\mathrm {Ann}}_{{\grdT}_{\grdQ}} f$ de $f$
dans $\grdT_{\grdQ}$ est $R^f :=\oplus_{g\neq f} R_g$, donc $J_f$ est
isog\`ene \`a $J_0 (M)_{\grdQ} /R^f J_0 (M)_{\grdQ}$. On a donc pour finir
l'isog\'enie :%
$$J_0 (N)_{\grdQ}  \to {\oplus}_{G_{\grdQ} .f} (J_f )^{\sigma_0 (N/M)} ,$$
o\`u la somme est prise sur les orbites sous l'action de Galois de toutes les
newforms en niveaux $M$ divisant $N$.

 L'accouplement bilin\'eaire :
$$\left\{ \begin{array}{rcl}
H_1 (X_0 (N )\, ;\, \grdQ )\times S_2 (N)_{\grdC}  & \to &\grdC \\
(c,f) &\mapsto &\langle c,f\rangle =\int_{c} f(z) dz
\end{array} \right. $$
est non d\'eg\'en\'er\'e, et pour lui les op\'erateurs de Hecke sont
auto-adjoints.

  On va d'abord prouver :
\begin{lemm}
Soit $f$ une newform en niveau $M$ divisant $N$ tel que $\langle e,f \rangle
\neq 0$. Alors ${\cal A}_e \cap R_f = \{ 0\} .$
\end{lemm}

{\bf Preuve du lemme.} Calculons l'effet de $B_D$ ($D$ entier quelconque) sur
l'accouplement :

$$ \langle e, B_D (f)\rangle = \int_{e} f (D.z)dz=i\int_{0}^{\infty} f
(iDz)dz =\frac{i}{D} \int_{0}^{\infty} f (iz)dz=\frac{1}{D} \langle
e,f\rangle .$$

   Avant de faire des calculs explicites, expliquons l'id\'ee de la preuve. On
a vu que $S_2 (N)_{\grdQ}$ \'etait un ${\grdT}_{\grdQ}$-module libre de rang
1 : comme tel, selon le th\'eor\`eme 3.5, il est isomorphe \`a un produit de
$\grdQ$-espaces vectoriels de type $K[X_1 ,\dots ,X_n ]/(X_1^{m_1} \dots
X_n^{m_n} )$, sur lesquels la restriction de ${\grdT}_{\grdQ}$ a la m\^eme
forme, et agit par simple multiplication. Si donc on prend un \'el\'ement
$t=\sum_{j} \lambda_j  \, X_{1}^{m_1^j} X_{2}^{m_2^j} \dots X_{m}^{m_n^j}$ de
la ``partie'' $K[X_1 ,\dots ,X_n ]/(X_1^{m_1} \dots X_n^{m_n} )$ de
${\grdT}_{\grdQ}~\cap~{\cal A}_e$, en le multipliant par des mon\^omes
judicieux de $K[X_1 ,\dots ,X_n ]/(X_1^{m_1} \dots X_n^{m_n} )$ (vu cette fois
comme le sous-espace des formes modulaires sur lequel $t$ agit), on aura que
l'accouplement de $e$ avec des mon\^omes de type $\lambda_j X_1^{\alpha_1}
\dots X_n^{\alpha_n}$ sera nul, pour tous les $\lambda_j$. Mais de tels
mon\^omes correspondront \`a des ``d\'ecalages'' par des op\'erateurs $B_d$ de
notre newform originelle $f$, d\'ecalages qui, comme le montre l'int\'egration
par partie ci-dessus, ne font que multiplier le produit $\langle e,f\rangle$
par
des constantes non nulles. On d\'eduira donc de la non nullit\'e de ce produit
$\langle e,f\rangle$ que tous les coefficients $\lambda_j$ sont nuls, et donc
$t$.

   \'Ecrivons-le maintenant pr\'ecis\'ement. Soit $t$ \'el\'ement de $R_f$. On
a les \'equivalen\-ces suivantes :%
$$(t.e=0)\Leftrightarrow (\forall g \in  S_2 (N)_{\grdC} , \langle t.e, g
\rangle =0)\Leftrightarrow (\forall g \in {\cal B}_f, \langle e, t.g\rangle
=0).$$ %
Supposons de plus $t$ dans ${\cal A}_e$. On d\'ecompose $E_f$ comme somme
directe des intersections des sous-espaces caract\'eristiques des op\'erateurs
$U_q$ pour $q$ divisant $N$. Consid\'erons la restriction $\tilde{t}$ de $t$
\`a l'un de ces sous-espaces, ${\tilde{E}}_f$. On \'ecrit $\tilde{t}$ comme un
polyn\^ome en les $U_q$, pour $q$ divisant $N$, \`a coefficients dans $\grdQ
[T_l ,...]$ ($l$ ne divisant pas $N$), c'est-\`a-dire \`a coefficients dans
(un corps isomorphe \`a) $K_f$. On peut ne consid\'erer que les $U_q$ qui
n'agissent pas diagonalement sur ${\tilde{E}}_f$. Si on note $\alpha_q$ la
valeur propre de $U_q$ sur ${\tilde{E}}_f$, alors $u_q :=(U_q -{\alpha}_q
{\mathrm {Id}}|_{{\tilde{E}}_f} )$ est nilpotent.

   On \'ecrit donc $\tilde{t}$ sous forme d'un polyn\^ome :
$$\tilde{t} =\sum_{j\geq 0}  \lambda_j  \, u_{q_1}^{\alpha_1^j}
u_{q_2}^{\alpha_2^j} \dots u_{q_m}^{\alpha_m^j} ,$$ %
o\`u les multi-indices $\alpha_k^j$ sont ordonn\'es par exemple de fa\c{c}on
lexicographique ($j<j'$ si $\alpha_s^j < \alpha_s^{j'} ,$ o\`u $s:={\mathrm
{inf}} \{ r$ tel que $\alpha_r^j \neq \alpha_r^{j'} \}$). %

   D'apr\`es la r\'eduction de Jordan de la partie pr\'ec\'edente, il existe
des polyn\^omes en $B_{q_j}$ : $P_j (B_{q_j} )$, et un \'el\'ement :
$$g_1 := [ P_1 (B_{q_1} )\otimes  P_2 (B_{q_2} )\otimes \cdots P_m (B_{q_m} )]
(f)\, \in S_2 (N)_{\grdC} \ ,$$ %
tel que $u_{q_j}^{\alpha_j^1} (P_j (B_{q_j} )(f))=: {\cal P}_j (B_{q_j} (f))
\neq 0$ et $u_{q_j}^{\alpha_j^1 +1} (P_j (B_{q_j} )(f))=0.$ Pour \^etre plus
explicite, ${\cal P}_j (B_{q_j} (f))$ a l'une des quatre formes suivantes : %
\begin{enumerate}
\item $(f)$ si on est sur l'espace caract\'eristique associ\'e \`a $0$ et
$q_j^2 |M$ ;
\item $(B_{q_j} (f) -a_{q_j} f)$ si on est sur l'espace
caract\'eristique associ\'e \`a $0$ et $q_j ||M$ ;
\item $(B_{q_j}^2 (f)-(a_{q_j} /{q_j} )B_{q_j} (f) +(1/q_j ) f)$ si on est sur
l'espace caract\'eristique associ\'e \`a $0$ et $q_j$ ne divise pas $M$ ;
\item $(-q_j B_{q_j} (f)+\alpha_{q_j} f)$ si $q_j$ ne divise pas $M$ et on est
sur le sous-espace Ker$(U^2_{q_j} -a_{q_j} U_{q_j} +q_j )=$Ker$(U_{q_j}
-\alpha_{q_j}
{\mathrm{Id}})^2$.
\end{enumerate}
Alors %
$$0=\langle e, t.g_1 \rangle =\lambda_1 \, \int_{0}^{i\infty} \prod_{j=1}^{m}
{\cal P}_j (B_{q_j} (f))\, dz\, =\lambda_1 \prod_{j=1}^{m} {\cal P}_j
(\frac{1}{q_j} ) \int_{0}^{i\infty} f\, dz$$ %
$$=\lambda_1 \prod_{j=1}^{m} {\cal P}_j (\frac{1}{q_j} )\langle e,f\rangle
\ ;$$   %
or aucun des facteurs ${\cal P}_j (\frac{1}{q_j} )$ n'est nul, puisque
respectivement de forme : %
\begin{enumerate}
\item 1,
\item $(1/q_j  -a_{q_j} )=(1/q_j \ \pm 1)$,
\item $(1/{q_j^2}-a_{q_j}/{q_j^2} \, +1/q_j )=(1/q_j^2 )(q_j \, +1\, -a_{q_j}
)$, ou encore :
\item $(\alpha_{q_j} -1)=(\pm \sqrt{q_j} -1).$
\end{enumerate}
(La seule non-nullit\'e qui ne soit pas triviale est la troisi\`eme : c'est la
borne de Weil ($|a_{q_j} |\leq 2\sqrt{q_j} )$ qui la montre.) Donc $\lambda_1$
est nul, et en recommen\c{c}ant la m\^eme op\'eration, par r\'ecurrence tous
les $\lambda_j$ sont nuls - donc $\tilde{t}$ l'est aussi. Ceci \'etant vrai
pour toutes les restrictions $\tilde{t}$ de $t$ aux intersections de
sous-espaces caract\'eristiques, $t$ est lui-m\^eme nul. $\Box$ \\
{\bf Fin de la preuve du th\'eor\`eme.}  On a montr\'e que si $\langle e,f
\rangle \neq 0$, ${\cal A}_{e,\grdQ} \cap R_f =0$, et r\'eciproquement il est
\'evident que si  $\langle e,f \rangle = 0$ alors ${\cal A}_{e,\grdQ} \cap R_f
=R_f$. Comme $L(f,1)= 2\pi \langle e,f\rangle$, on peut donc \'ecrire : %
$${\cal A}_{e,\grdQ} = \bigoplus_{G_{\grdQ} f /L(f,1)= 0} R_f \ ,$$ %
et on a des isog\'enies : %
$$J^e_0 (\grdQ )=(J_0 (N)/{\cal A}_e J_0 (N))(\grdQ ) \to \prod_{R_f
\not\subset {\cal A}_{e,\grdQ}} (J_0 (N)/({\grdT}_{\grdQ} /R_f )J_0 (N)) (\grdQ
)
\rightarrow$$ %
$$\rightarrow \prod_{G_{\grdQ} f /L(f,1)\neq 0} J_f (\grdQ ).$$ %
Selon le th\'eor\`eme de Kolyvagin-Logachev, les $J_f (\grdQ )$ ci-dessus
sont finies (\`a cause justement de la non-nullit\'e en 1 de la fonction
$L(f,s)$, voir le th\'eor\`eme 0.3 de \cite{Kolyvagin}, compl\'et\'e par des
r\'esultats de \cite{Murty-Murty} ou \cite{Bump-Friedberg-Hoffstein}). Et donc
le quotient d'enroulement est bien de rang nul. $\Box$
\section{Espaces tangents et cotangents.}
\subsection{Espace tangent \`a $J_0 (N)_{/\grdZ [1/N]}$ en z\'ero.}
Soit $N$ un entier, fix\'e pour toute cette partie. On montre dans cette
sous-section :
\begin{theo}
L'espace tangent \`a $J_0 (N)_{/\grdZ [1/N]}$ en z\'ero, ${\mathrm {Tan}}_0
(J_0 (N)_{/\grdZ [1/N]})$, est un ${\grdT}_{\grdZ} \otimes \grdZ [1/N]$-module
libre de rang $1$, de base ${\frac{d}{dq} |}_0$.
\end{theo}
On va donner quatre lemmes pr\'eliminaires, desquels le th\'eor\`eme
d\'ecoulera naturellement. Dans ce qui suit, pour all\'eger les notations, on
notera $M$ le module ${\mathrm {Tan}}_0 (J_0 (N)_{/\grdZ [1/N]} )$, et $R$
l'anneau ${\grdT}_{\grdZ} \otimes \grdZ [1/N]$, qui est un $\grdZ
[1/N]$-module libre de type fini (${\grdT}_{\grdZ}$ l'est sur ${\grdZ}$,
puisqu'on peut le voir comme un sous-module des endomorphismes de $J_0
(N)_{/\grdZ}$). %
\begin{lemm}
Le $R$-module $M$ est fid\`ele.
\end{lemm}
{\bf Preuve.} Comme ni $R$ ni $M$ n'ont de $\grdZ [1/N]$-torsion, on peut voir
le lemme en \'etendant les scalaires \`a $\grdQ$, et en remarquant qu'alors,
le dual de $M\otimes \grdQ$ (comme $\grdQ$-espace vectoriel), qui est
${\mathrm {Cot}}_0 (J_0(N)_{\grdQ} )$, est bien $\grdT \otimes
\grdQ$-fid\`ele, comme on le voit dans la partie pr\'ec\'edente. $\Box$
\begin{lemm}
Pour tout id\'eal maximal ${\cal M}$ de $R$, on a $M/{\cal M}M\neq 0$.
\end{lemm}
{\bf Preuve.} Supposons qu'il existe un id\'eal maximal ${\cal M}$ tel que
$M/{\cal M}M=0$. Puisque $M$ est un $R$-module fini, le faisceau associ\'e
$\tilde{M}$ sur le sch\'ema affine noeth\'erien $X:=$Spec($R$) est coh\'erent
; notons $x$ le point correspondant \`a ${\cal M}$. On a :
$$M/{\cal M}M=M\otimes_R R/{\cal M}=M\otimes_R R_{\cal M} /{\cal M}R_{\cal M}=
M_{\cal M} /{\cal M}M_{\cal M} ,$$
donc le lemme de Nakayama dit que $M_{\cal M} =0$. Si $g$ est un point
g\'en\'erique de $X$ (correspondant \`a l'id\'eal ${\cal P}$ de $R$) qui se
sp\'ecialise en $x$, on a $\tilde{M}_g =0$ (comme localisation de
$\tilde{M}_x$). Puisque $R$ est sans $\grdZ$-torsion, $g$ est au-dessus du
point g\'en\'erique de $\grdZ [1/N]$ ; et puisque c'est un $\grdZ
[1/N]$-module de type fini, Spec($R\otimes_{\grdZ [1/2N]} \grdQ$) est de
dimension nulle. On a donc $R\otimes \grdQ \simeq \prod_{\cal Q} R_{\cal Q}$,
o\`u le produit (fini) est pris sur les id\'eaux minimaux ${\cal Q}$ de
$R$ : ce qui veut dire que $R_{\cal P}$ est isomorphe \`a un facteur de
$R\otimes_{\grdZ [1/N]} \grdQ$. Mais alors la fid\'elit\'e de $M$ comme
$R$-module contredit le fait que $\tilde{M}_g =0$. $\Box$
\begin{lemm}
Pour tout id\'eal maximal ${\cal M}$ de $R$, l'\'element ${\frac{d}{dq} |}_0$
de $M$ est d'image non nulle dans le quotient $M/{\cal M}M$.
\end{lemm}
{\bf Preuve.} Soit ${\cal M}$ un id\'eal maximal de $R$. De la finitude comme
$\grdZ$-module de ${\grdT}_{\grdZ}$ on d\'eduit d'abord que le sous corps
premier de $F:=R/{\cal M}$ est n\'ec\'essairement fini, ${\grdF}_l$, et que
$F$ en est une extension finie. Posons $\overline{R} :=R\otimes_{\grdZ}
{\overline{\grdF}}_l$, et $\overline{M} :=M\otimes_{\grdZ}
{\overline{\grdF}}_l$. On a la suite exacte : ${\cal M} \otimes_{\grdZ}
{\overline{\grdF}}_l \to \overline{R} \to F\otimes_{\grdZ}
{\overline{\grdF}}_l \to 0.$ Choisissons un plongement $i$ : $F
\hookrightarrow {\overline{\grdF}}_l$, et d\'eduisons-en un morphisme de
${\overline{\grdF}}_l$-alg\`ebre : %
$$i\otimes 1\, :\, F\otimes_{\grdF_l} {\overline{\grdF}}_l \to
{\overline{\grdF}}_l .$$
Notons $\overline{{\cal M}}$ son noyau : on a la suite exacte $0\to
\overline{{\cal M}} \to \overline{R} \to {\overline{\grdF}}_l \to 0$. On
v\'erifie alors qu'on a :
$$\left.\frac{d}{dq}\right|_0 \mapsto (M/{\cal M}M)\hookrightarrow (M
/{\cal M}M) \otimes_{F} {\overline{\grdF}}_l \simeq \overline{M} /
\overline{{\cal M}} \, \overline{M} .$$
Si $E$ est un ${\overline{\grdF}}_l$-espace vectoriel, notons ${E}^{\vee}$ son
dual, et d\'esignons par ${\overline{M}}^{\vee} [\overline{{\cal M}} ]$ le
sous-$\overline{R}$-module de ${\overline{M}}^{\vee}$ tu\'e par
$\overline{{\cal M}}$. Avec ces notations,
$$(\overline{M} /\overline{{\cal M}} \, \overline{M} )^{\vee} \simeq
\overline{M}^{\vee} [\overline{{\cal M}} ]$$
en tant que $\overline{R}$-modules. (En effet, les formes
${\overline{\grdF}}_l$-lin\'eaires sur l'espace vectoriel quotient
$(\overline{M} /\overline{{\cal M}} \, \overline{M} )$ s'identifient avec
celles sur $\overline{M}$ dont le noyau contient $\overline{{\cal M}} \,
\overline{M}$, et l'action de $\overline{R}$-module que l'on consid\`ere est :
$(t.f)(\ )=f(t.\ ),$ ($t\in \overline{R}$, $f\in {\overline{M}}^{\vee}$).) On
a donc enfin :
$$\overline{M}^{\vee} [\overline{{\cal M}} ]=H^0 (X_0 (N
)_{{\overline{\grdF}}_l}, {\Omega}^1 )[\overline{{\cal M}} ]\simeq S_2 (N
)_{\overline {\grdF}_l} [\overline{{\cal M}} ].$$
On a montr\'e dans le lemme pr\'ec\'edent que $M/{\cal M}M$ \'etait non nul :
soit donc $f$ une forme parabolique non nulle appartenant \`a $S_2
(N)_{\overline {\grdF}_l} [\overline{{\cal M}} ]$. Le corollaire III, 12.9 de
\cite{Hartshorne} assure que $H^0 (X_0 (N)_{{\overline{\grdF}}_l} ,{\Omega}^1
)=H^0 (X_0 (N)_{/\grdZ [1/N]} ,{\Omega}^1 )\otimes {\overline{\grdF}}_l$, ce
qui implique que, quitte \`a prendre un rel\`evement de $f$, on peut supposer
qu'on est en caract\'eristique nulle, et m\^eme dans $S_2 (N)_{\grdC}$ : pour
voir que les coefficients de $f$ v\'erifient la relation : %
$$a_1 (T_n .f)=a_n (f)$$ %
m\^eme si $f$ n'est pas une newform (cela r\'esulte par exemple de la formule
(3.5.12) de \cite{Shimura}). On en d\'eduit que pour tout entier $n$,
$\frac{d}{dq} |_0 (T_n .f) =a_n (f)$. Mais puisque $X_0 (N)_{{\overline
\grdF}_l}$ est int\`egre et lisse, ses formes diff\'erentielles sont
d\'efinies par leur d\'eveloppement de Fourier en l'infini : donc
$\frac{d}{dq} |_0$ ne peut \^etre d'image nulle. $\Box$ %
\begin{lemm}
Pour tout id\'eal maximal ${\cal M}$ de $R$, $M/{\cal M}M$ est un
$R/{\cal M}$-module libre de rang $1$, de base l'image de $\frac{d}{dq} |_0$.
\end{lemm}
{\bf Preuve.} Les r\'esultats pr\'ec\'edents montrent qu'il suffit maintenant
de prouver que ${\mathrm{dim}}_{R/{\cal M}} (M/{\cal M}M)\leq 1$. Si $(l)$ est
l'id\'eal premier de $\grdZ$ au-dessus duquel ${\cal M}$ se trouve,
choisissons comme dans la preuve du pr\'ec\'edent lemme un plongement de
$F=R/{\cal M}$ dans ${\overline{\grdF}}_l$. Toujours selon la preuve
pr\'ec\'edente et avec les m\^emes notations, on a les \'equivalences : %
$$((M/{\cal M}M)\ {\mathrm {est\  un\  }}F{\mathrm {-espace\ vectoriel\ de\
dimension\ }}1)$$ %
\centerline{$\Updownarrow$} %
$$((M/{\cal M}M)\otimes_{F} {\overline {\grdF}}_l \ {\mathrm {est\ un\ }}
{\overline \grdF}_l {\mathrm {-espace\ vectoriel\ de\ dimension\ }} 1)$$ %
\centerline{$\Updownarrow$} %
$$((\overline{M} /\overline{{\cal M}} \, \overline{M})\ {\mathrm {est\ un\ }}
{\overline \grdF}_l {\mathrm {-espace\ vectoriel\ de\ dimension\ }} 1)$$ %
\centerline{$\Updownarrow$} %
$$({\overline{M}}^{\vee} [\overline{{\cal M}} ]\ {\mathrm {est\ un\ }}
{\overline \grdF}_l {\mathrm {-espace\ vectoriel\ de\ dimension\ }}1).$$ %
Pour tout entier $n$, soit $a_n$ l'image dans $\overline{R} /\overline{{\cal
M}} \simeq {\overline \grdF}_l$ de l'op\'erateur de Hecke $T_n$. Soit $f$ un
\'el\'ement non nul de $H^0 (X_0 (N)_{{\overline \grdF}_l} , {\Omega}^1
)[\overline{{\cal M}} ]$ ; $f$ est une forme propre pour les op\'erateurs de
Hecke : pour tout $n$, $T_n (f)=a_n f$. \`A une constante multiplicative non
nulle pr\`es, le d\'eveloppement de Fourier de $f$ en l'infini est donc :
$q+a_2 q^2 +a_3 q^3 +\cdots$. De m\^eme que plus haut, on d\'eduit du fait que
$X_0 (N)_{{\overline \grdF}_l}$ soit int\`egre et lisse que toute forme
diff\'erentielle sur $X_0 (N)_{{\overline \grdF}_l}$ est d\'efinie par son
$q$-d\'eveloppement en l'infini ; donc $f$ engendre l'espace vectoriel
$S_2(N)_{\overline {\grdF}_l}.\ \Box$ \\ %
{\bf Preuve du th\'eor\`eme.} Montrons que le module $M':=M/(R.(\frac{d}{dq} )
|_0 )$ est nul. D'apr\`es le lemme pr\'ec\'edent, pour tout id\'eal maximal
${\cal M}$ de $R$, %
$$M'/{\cal M}M'=M'\otimes_R R_{\cal M} /{\cal M}R_{\cal M}=0\ ;$$ %
donc par le lemme de Nakayama, $M'_{\cal M}$ est nul. On peut alors appliquer
le raisonnement de la preuve du lemme 4.4, et en conclure que le faisceau
$\tilde{M'}$ sur Spec($R$) est de fibre nulle en chaque point, donc nul.
Enfin, $M$ est un $R$-module libre puisque fid\`ele. $\Box$ %

  Notons au passage un corollaire de ce dernier th\'eor\`eme (voir aussi
3.7) :
\begin{prop}
L'espace des formes paraboliques $S_2 (N)_{\grdQ}$ est un
$\grdT_{\grdQ}$~-~module libre de rang 1.
\end{prop}
{\bf Preuve.} Une application du corollaire 3.6 du chapitre pr\'ec\'edent
donne que la $\grdQ$-alg\`ebre de Hecke $\grdT_{\grdQ}$ (pour ${\Gamma}_0
(N)$, en poids 2) est de Gorenstein. De plus, on vient de montrer que
${\mathrm {Tan}}_0 (J_0 (N)_{\grdQ} )$ \'etait comme $\grdT_{\grdQ}$-module
isomorphe \`a $\grdT_{\grdQ}$ lui-m\^eme ; et comme ${\mathrm {Cot}}_0 (J_0
(N)_{\grdQ} )={\mathrm {Tan}}_0 (J_0 (N)_{\grdQ} )^{\vee}$, ${\mathrm {Cot}}_0
(J_0 (N)_{\grdQ} )$ est isomorphe \`a $\grdT_{\grdQ}$ soi-m\^eme comme
$\grdT_{\grdQ}$-module. $\Box$ %
\subsection{Espace tangent au quotient de la jacobienne.}
Soit $I$ un id\'eal de $\grdT_{\grdZ}$, satur\'e ({\em i.e.} tel que $\grdT/I$
soit sans $\grdZ$-torsion). On d\'eduit de la jacobienne $J_0 (N)$ sur
${\grdZ [1/2N]}$ les sch\'emas ab\'eliens $J^I :=I.J_0 (N)$ et $J_I :=J_0 (N)
/I.J_0 (N).$ Le but de cette section est de d\'emontrer le th\'eor\`eme
suivant, voisin et corollaire de la celui de la section pr\'ec\'edente :%
\begin{theo}
L'espace tangent \`a $J_{I,\grdZ [1/2N]}$ en z\'ero, ${\mathrm {Tan}}_0
(J_{I,\grdZ [1/2N]})$, est un $\grdT /I \otimes \grdZ [1/2N]$-module libre
de rang 1, de base l'image de $\frac{d}{dq} |_0$.
\end{theo}
On commence par travailler sur $\grdQ$ :
\begin{prop}
On a ${\mathrm {Tan}}_0 (J^I_{\grdQ} )=I.{\mathrm {Tan}}_0 (J_0 (N)_{\grdQ}
)$.
\end{prop}
{\bf Preuve de la proposition.} Soit $(i_1 ,i_2 ,\dots ,i_n )$ un syst\`eme de
g\'en\'erateurs de $I$ ($\grdT$ est noeth\'erien puisque $\grdZ$-module de
type fini) ; $J^I_{\grdQ}$ est donc l'image du morphisme :
$$\oplus_{1\leq j\leq n} J_0 (N)_\grdQ \stackrel{\phi}{\rightarrow}
J_0 (N)_\grdQ$$
d\'efini par la multiplication par ce syst\`eme. On en d\'eduit sur les
espaces tangents :%
$$\oplus_{1\leq j\leq n} {\mathrm {Tan}}_0 (J_0 (N)_{\grdQ} ) \stackrel{
\tilde{\phi}}{\rightarrow} {\mathrm {Tan}}_0 (J_0 (N)_{\grdQ} ).$$
La proposition d\'ecoule alors du lemme :%
\begin{lemm}
Pour tout morphisme $A
\stackrel{\psi}{\rightarrow} B$ de vari\'et\'es ab\'eliennes sur $\grdQ$, on
a ${\mathrm {Tan}}_0 (\psi (A)) =\tilde{\psi} ({\mathrm {Tan}}_0 (A)).$
\end{lemm}
{\bf Preuve du lemme.}  Le foncteur $A \mapsto {\mathrm {Tan}}_0 (A)$ de la
cat\'egorie des groupes alg\'ebriques commutatifs sur $\grdQ$ dans celle des
$\grdQ$-espaces vectoriels, est exact \`a gauche. Mais il est aussi exact \`a
droite : de la suite exacte $0\to A\to B\to C\to 0$ se d\'eduit celle sur les
tangents : $0\to {\mathrm {Tan}}_0 (A)\to {\mathrm {Tan}}_0 (B)\to {\mathrm
{Tan}}_0 (C)$ ; si les dimensions des vari\'et\'es sont $a,\ b$ et $c$
respectivement, on a $a+c=b$, et comme les dimensions des espaces tangents
associ\'es sont les m\^emes (lissit\'e des groupes alg\'ebriques sur $\grdQ$),
la derni\`ere fl\`eche de la derni\`ere suite est surjective.

  Si donc on a un morphisme $\psi$ :%
$$0\to {\mathrm {ker}}(\psi )\to A\stackrel{\psi}{\rightarrow} \psi (A)\to
0,$$ %
on a bien %
$${\mathrm {Tan}}_0 (\psi (A))={\mathrm {Tan}}_0 (A)/{\mathrm {Tan}}_0
({\mathrm {ker}}(\psi ))= \tilde{\psi} ({\mathrm {Tan}}_0 (A) ).\ \Box$$ %

  Consid\'erons maintenant le probl\`eme sur $\grdZ [1/2N]$.\\ %
{\bf Preuve du th\'eor\`eme.}  De la suite exacte de vari\'et\'es ab\'eliennes
sur $\grdQ$ :%
$$0\to J^I_{\grdQ} \to J_0 (N)_{\grdQ} \to J_{I,\grdQ} \to 0,$$ %
on d\'eduit par propri\'et\'e universelle des mod\`eles de N\'eron un complexe
de sch\'emas ab\'eliens : %
$$0\to J^I_{\grdZ [1/2N] } \to J_0 (N)_{\grdZ [1/2N]} \to J_{I,\grdZ [1/2N]}
\to 0,$$ %
qui est en r\'ealit\'e une {\em suite exacte}, d'apr\`es un r\'esultat de
Raynaud (voir section 1, proposition 1.2 de \cite{rational}). Ce r\'esultat
assure encore que la suite :%
$$0\to {\mathrm {Tan}}_0 (J^I_{\grdZ [1/2N] }) \to  {\mathrm {Tan}}_0
(J_0 (N)_{\grdZ [1/2N]} )\to {\mathrm {Tan}}_0 (J_{I,\grdZ [1/2N]} )\to 0$$ %
est exacte elle aussi (m\^eme r\'ef\'erence, corollaire 1.1). Or on voit par
la section pr\'ec\'edente que ${\mathrm {Tan}}_0 (J_0 (N)_{\grdZ [1/2N]} )$
est un $\grdT \otimes \grdZ [1/2N]$-module libre de rang 1, de base
$\frac{d}{dq} |_0$ ; il existe donc un id\'eal $I'$ de $\grdT [1/2N]$ tel que
la suite pr\'ec\'edente de $\grdT [1/2N]$-modules soit isomorphe \`a :%
$$0\to I'\to \grdT \otimes \grdZ [1/2N] \to \grdT \otimes \grdZ [1/2N] /I'\to
0.$$ %
Mais puisque sur $\grdQ$ on a $I_{\grdQ} =I'_{\grdQ}$, que $I$ est par
hypoth\`ese satur\'e, et que $I'$ l'est aussi (puisque ${\mathrm {Tan}}_0
(J_{I,\grdZ [1/2N]} )\simeq \grdT [1/2N]/I'$ est sans torsion), on a bien
$I'=I_{\grdZ [1/2N]}$. Donc %
$${\mathrm {Tan}}_0 (J_{I,\grdZ [1/2N]} )\simeq {(\grdT /I)}_{\grdZ
[1/2N]} .\ \Box $$ %
\subsection{Preuve de la proposition ``crit\`ere de Kamienny''.}
On rappelle la situation dans laquelle on s'est mis avec notre probl\`eme
initial de borne pour la torsion de courbes elliptiques (voir la section 1,
dont on reprend les hypoth\`eses et notations). Rappelons le morphisme naturel
$f_d$ : ${X_0 (p^n )}^{(d)}_{\rm lisse} \rightarrow J_0^e $,
normalis\'e par ${\infty }^{(d)} \mapsto 0$. On a d\'efini en (1.2) un point
$j'^{(d)}$, \`a valeur dans $\grdZ$, de $X_ 0 (p^n )^{(d)}_{\rm lisse}$,
qui croise $\infty^{(d)}$ sur la fibre en $l$. On va d'abord montrer que ceci
implique que $f_d$ n'est pas une immersion formelle au point
${\infty}^{(d)}_{\grdF_l}$, en suivant l'expos\'e de \cite{Oest}.

  Soit $S$ un sch\'ema quelconque. Si $\phi$ : $X\to Y$ est un morphisme de
$S$-sch\'emas noeth\'eriens, on dit que c'est une {\em immersion formelle} en
un point $x$ de $X$ si l'application ${\hat{\cal O}}_{Y,\phi (x)} \to
{\hat{\cal O}}_{X,x}$ d\'eduite de $\phi$ entre compl\'et\'es d'anneaux locaux
est surjective. Si $s$ est l'image dans $S$ de $x$ (et $y$), on montre que
ceci \'equivaut \`a ce que la restriction de $\phi$ aux fibres en $s$ soit une
immersion formelle en $x$. On d\'emontre ensuite dans \cite{Oest} %
(lemme 5.1) : %
\begin{lemm}
Supposons que $X$ soit s\'epar\'e, que $\phi$ : $X\to Y$ soit une immersion
formelle en $x$. Supposons qu'il existe un sch\'ema int\`egre noeth\'erien $T$
et deux points $p_1$ et $p_2$ de $X$ \`a valeur dans $T$, tel qu'en un point
$t$ de $T$ on ait $x=p_1 (t)=p_2 (t)$. Si de plus on a $\phi \circ p_1 =\phi
\circ p_2$, alors $p_1 =p_2$.
\end{lemm}
Ce qui nous permet de montrer :
\begin{theo}
Le morphisme $f_d$ n'est pas une immersion formelle en ${\infty}^{(d)}_{
\grdF_l}$.
\end{theo}
{\bf Preuve.} Puisque $j'^{(d)}$ croise $\infty^{(d)}$ sur la fibre en $l$, le
point $f_d (j'^{(d)} )$ de $J_0^e ({\grdZ}_l )$ se r\'eduit en $0$ au-dessus
de $l$. Ce point est de torsion, car $J_0^e (\grdQ )$ l'est. Puisque $l>2$, le
``lemme de sp\'ecialisation'' (voir par exemple \cite{Bas}, lemme 3.3) assure
que $f_d (j'^{(d)} )$ est nul dans tout $J_0^e ({\grdZ}_{(l)} )$. Mais alors,
on peut appliquer le lemme pr\'ec\'edent et en conclure que si $f_d$ \'etait
une immersion formelle en ${\infty}^{(d)}_{\grdF_l}$, on aurait $j'^{(d)} =
{\infty}^{(d)}$ - ce qui contredit \'evidemment l'interpr\'etation modulaire
de ces points. $\Box$

   La courbe $X_0 (N)$ a en %
l'infini la coordonn\'ee formelle $q$ ; donc $X_0 (N)^d$ a les coordon\-n\'ees
formelles $q_1,\dots ,\, q_d$ au point $(\infty ,\dots ,\, \infty )$, et
les fonctions sym\'etriques \'el\'ementaires $\sigma_1 =q_1 + \cdots +q_d ,\,
\dots ,\ \sigma_d =q_1 \cdots q_d$ sont des coordon\-n\'ees formelles au point
${\infty}^{(d)}$ de $X_0 (N)^{(d)}$, donc $d\sigma_1 ,\dots ,\, d\sigma_d$
forment une base de ${\mathrm {Cot}}_{\infty^{(d)}} (X_0 (N)^{(d)} )$. On a
alors le lemme : %
\begin{lemm}
Soit $\omega$ un \'el\'ement de ${\mathrm {Cot}}_0 (J_0^e )$. Alors $f_1^*
(\omega )$ est une forme diff\'e\-ren\-tielle sur $X_0 (N)_{/\grdZ [1/N]}$.
Notons
$(\sum_{n\geq 1} a_n q^n )(dq/q)$ son d\'eveloppement de Fourier \`a l'infini.
On a :
$$f_d^* (\omega )=a_1 d\sigma_1 -a_2 d\sigma_2 + \cdots +(-1)^{d-1} d\sigma_d
\ \in {\mathrm {Cot}}_{\infty^{(d)}} (X_0 (N)^{(d)}_{/\grdZ [1/N]} ).$$
\end{lemm}
Pour la {\bf preuve}, voir \cite{Bas}, lemme 4.2. On en d\'eduit :
\begin{cor}
L'application tangente \`a $f_d$ en $\infty^{(d)}$ envoie $\frac{d}{d\sigma_1}
,\, \dots ,\, \frac{d}{d\sigma_d}$ sur $T_1 (\frac{d}{dq} |_0 ) , -T_2
(\frac{d}{dq} |_0 ),\, \dots ,\, (-1)^{d-1} T_d (\frac{d}{dq} |_0 )$
respectivement.
\end{cor}
{\bf Preuve.} La relation $a_1 (T_n f)=a_n (f)$ pour toute forme modulaire
(voir la preuve du lemme 4.5) et le lemme pr\'ec\'edent suffisent \`a
conclure. $\Box$

  Et on arrive au th\'eor\`eme qu'on veut finalement montrer :
\begin{theo}
 On a \'{e}quivalence entre :
\begin{enumerate}
\item $f_d$ est une immersion formelle en ${\infty}_{{\grdF}_l}^{(d)}$, et
\item $T_1 e,...,T_d e$ sont $\grdF _l$-lin\'{e}airement ind\'{e}pendants dans
$\grdT e/l\grdT e$.
\end{enumerate}
De plus, ces deux conditions sont satisfaites si l'est :
\begin{enumerate}
\item[3.] $T_1 \{ 0,\infty \} ,...,T_{d.s} \{ 0,\infty \}$ sont
$\grdF_l$-lin\'{e}airement ind\'{e}pendants dans l'espace vectoriel
$H_1 (X_0 (p^n ),\, {\mathrm {pointes }}\, ;\, \grdZ )\otimes \grdF_l$ (o\`u
$s$ d\'{e}signe le plus petit nombre premier diff\'erent de $p$).
\end{enumerate}
\end{theo}
Avant de passer \`a la d\'emonstration, un lemme pr\'eliminaire :%
\begin{lemm}
Soit $N$ un entier positif strict, et $r$ un second, premier au premier.
Alors, sur la courbe modulaire $X_0 (N)_{\grdC}$, on a $T_r .0={\sigma}_1 (r).
0$ et $T_r .\infty ={\sigma}_1 (r).\infty ,$ o\`u comme d'habitude ${\sigma}_1
(r)$ est la somme des diviseurs de $r$.
\end{lemm}
{\bf Preuve du lemme.} On se place ``en projectif'' : on identifie $X_0
(N)_{\grdC}$ avec $(\grdH \cup {\grdP}^1 (\grdQ ))/{\Gamma}_0 (N)$ et on
d\'esigne par $({x\atop y})$ l'\'el\'ement $x/y$ de ${\grdP}^1 (\grdQ )$. En
particulier, les points $0$ et $\infty$ seront \'ecrits $({0\atop 1})$ et
$({1\atop 0})$ respectivement. On peut voir les op\'erateurs de Hecke comme
des correspondances sur  $X_0 (N)_{\grdC}$ (voir la section 5), qui
v\'erifient : %
$$T_r  .0=\sum_{\begin{array}{c}
1\leq \delta |r\\
0\leq \beta <\delta
\end{array} }
\left( \begin{array}{cc}
\frac{r}{\delta} & -\beta \\
0 & \delta
\end{array}  \right) .\left( \begin{array}{c} 0 \\ 1 \end{array} \right)
=\sum_{\begin{array}{c}
1\leq \delta |r\\
0\leq \beta <\delta
\end{array} } \left( \begin{array}{r} -\beta  \\ \delta \end{array} \right).$$
Or, chacun de ces $({-\beta \atop \delta})$ est conjugu\'e par un \'el\'ement
de ${\Gamma}_0 (N)$ \`a $({0\atop 1})$ : en effet, quitte \`a remplacer
$\beta$ et $\delta$ par leur quotient par leur p.g.c.d., on peut les supposer
premiers entre eux, et choisir des entiers $a$ et $b$ tels que $a\delta -bN
\beta =1$ ; alors on a :%
$$\left( \begin{array}{cc}
a & -\beta \\
bN & \delta
\end{array} \right).\left( \begin{array}{c} 0 \\ 1 \end{array} \right) =\left(
\begin{array}{c} -\beta  \\ \delta \end{array} \right) ,$$ %
et la matrice de gauche est bien de ${\Gamma}_0 (N)$.

  De la m\^eme fa\c{c}on, %
$$T_r  .\infty =\sum_{\begin{array}{c}
1\leq \delta |r\\
0\leq \beta <\delta
\end{array} }
\left( \begin{array}{cc}
\frac{r}{\delta} & -\beta \\
0 & \delta
\end{array} \right) .\left( \begin{array}{c} 1 \\ 0 \end{array} \right)
=\sum_{\begin{array}{c}
1\leq \delta |r\\
0\leq \beta <\delta
\end{array} } \left( \begin{array}{c} \frac{r}{\delta}  \\ 0 \end{array}
\right) = \sigma_1 (r) .\infty \ \Box$$ %
{\bf Preuve du th\'eor\`eme.}\ D'abord l'\'equivalence des deux premiers
points. %
La propri\'et\'e pour un morphisme d'\^etre une immersion formelle en un
point est, formulation duale de celle sur les espaces cotangents, que
l'application induite par ce morphisme sur les tangents correspondants est une
injection. Le corollaire 4.16 dit donc que $f_d$ est une immersion formelle en
${\infty}^{(d)}_{\grdF_l}$ si et seulement si les vecteurs $T_i
({\frac{d}{dq}}|_0 )$, $1\leq i\leq d$ sont lin\'eairement ind\'ependants %
dans ${\mathrm {Tan}}_0 (J_0^e %
({\overline{\grdF}}_l ))$. Mais le th\'eor\`eme 4.9 assure que cet espace
tangent est un $\grdT /{\cal A}_e \otimes_{\grdZ} {\overline{\grdF}}_l$-module
libre de rang 1 de base $\frac{d}{dq} |_0$, d'o\`u l'\'equivalence des
propri\'et\'es {\em 1.} et {\em 2.} du th\'eor\`eme.

   Prouvons la derni\`ere implication. On rappelle qu'il y a un isomorphisme
de $\grdR$-espaces vectoriels : %
$$\left\{ \begin{array}{rcl}
H_1 (X_0 (p^n ) \, ;\, \grdZ )\otimes \grdR & \to & {\rm Hom}_{\grdC }
\left( H^0 (X_0 (p^n )\, ;\, {\Omega }^1 ),\ \grdC \right) \\
\gamma \otimes 1 & \mapsto & \left( \omega \mapsto \int_{\gamma } \omega
\right) .
\end{array} \right.$$ %
L'int\'egration permet d'en d\'eduire un morphisme de $H_1 (X_0 (p^n ),\,
{\rm pointes}\, ;\, \grdZ )$ dans $H_1 (X_0 (p^n ) \, ;\, \grdZ )\otimes
\grdR$, et comme on l'a dit, l'\'el\'ement d'enroulement $e$, qui est l'image
par ce morphisme du symbole modulaire $\{ 0,\infty \}$, est \'el\'ement de
$H_1 (X_0 (p^n ) \, ;\, \grdQ )$.

   Maintenant, c'est avec les symboles modulaires, vivant dans l'homologie
relative aux pointes, qu'on sait travailler ; mais c'est dans l'homologie
absolue qu'il nous faut une ind\'ependance $\grdF_l$-lin\'eaire. Pour s'y
ramener, on fait agir l'op\'erateur $I_s :=T_s -{\sigma}_1 (s)$ sur $\{ 0,
\infty\}$. En effet, puisque $s$ est premier au niveau, le lemme pr\'ec\'edent
dit que $T_s$ envoie les pointes $0$ et $\infty$ sur $(s+1)$-fois
elles-m\^emes, donc $I_s$ pousse bien $\{ 0,\infty \}$ dans $H_1 (X_0 (p^n )
\, ;\, \grdZ )$ vu comme sous-module de $H_1 (X_0 (p^n ),\, {\rm pointes}
\, ;\, \grdZ )$. Consid\'erons aussi $H_1 (X_0 (p^n )\, ;\, \grdZ )$ comme un
sous-module de $H_1 (X_0 (p^n )\, ;\, \grdQ )$. %

   Supposons que la propri\'et\'e {\em 2.} du th\'eor\`eme ne soit pas
v\'erifi\'ee. Il existe donc une relation de d\'ependance : %
$$\overline{\lambda_1} T_1 e+\cdots +\overline{\lambda_k} T_k e =0 \in \grdT
e/l\grdT e,$$ %
pour un $k\leq d$ tel que $\overline{\lambda_k}$ soit non nul. On la rel\`eve
en %
$$\lambda_1 T_1 e+\cdots +\lambda_k T_k e = l\, x \in l\grdT e \subseteq H_1
(X_0 (p^n )\, ;\, \grdQ ),$$ %
on la multiplie par $I_s$ pour obtenir %
$$I_s (\sum_{i=1}^{k} {\lambda}_i T_i e)= \sum_{i=1}^{k} {\lambda}_i T_i (I_s
e)={\lambda}_k T_{sk} \, e+\sum_{i=1}^{sk-1} {\mu}_i T_i e =l\, I_s x \in l
\grdT e.$$ %
(En effet, la d\'efinition des op\'erateurs de Hecke (par exemple avec la
s\'erie formelle) donne que $T_s .T_n =T_{s.n} +$ (combinaison lin\'eaire de
termes d'indice plus petit).) Les formes lin\'eaires que d\'efinissent $I_s
\{ 0, \infty \}$ et $I_s e$ par l'int\'egration sont les m\^emes, et puisqu'on
peut voir ces deux \'el\'ements comme appartenants \`a $H_1 (X_0 (p^n ) \, ;\,
\grdR )$, ils sont en fait \'egaux. On en d\'eduit une relation de
d\'ependance lin\'eaire : %
$${\lambda}_k T_{sk} \, \{ 0,\infty \} +(\sum_{i=1}^{sk-1} {\mu}_i T_i \{ 0,
\infty \} )=l\, I_s x \in lH_1 (X_0 (p^n ),\, {\rm pointes} \, ;\, \grdZ
),$$ %
et en consid\'erant son image dans $H_1 (X_0 (p^n ),\, {\rm pointes} \, ;
\, \grdZ )\otimes {\grdF}_l$, une contradiction avec la propri\'et\'e {\em 3.}
du th\'eor\`eme. $\Box$
\section{Lemme combinatoire.}
Le but de cette partie est de prouver la proposition-cl\'e de la
d\'emonstration g\'en\'erale :\\
{\bf Proposition 1.7}\ {\it Soit $p$ un nombre premier. Posons $C=\sqrt{65}$
si $p$ est diff\'erent de $2$, et $C=\sqrt{129}$ si $p$ est $2$. Notons $s$ le
plus petit nombre premier diff\'erent de $p$. Si $p^n>C^2.(sd)^6$, alors les
$T_i \{ 0,\infty \} ,\ 1\leq i\leq sd$ sont $\grdF$-lin\'{e}airement
ind\'{e}pendants dans le $\grdF$-espace vectoriel $H_1 (X_0 (p^n ),\,
{\mathrm {pointes }}\, ;\, \grdZ )\otimes \grdF$ pour tout corps $\grdF$.}
\subsection{Notations et rappels.}
On identifie $\Gamma_0(p^n) \backslash SL_2 (\grdZ)$ \`a ${\grdP}^1 \left(
{\grdZ}/{p^n}{\grdZ} \right)$ avec : %
$$ {\Gamma }_0 ({p^n})\left( \begin{array}{cc}
a & b \\
c & d
\end{array}
\right)
\mapsto (\overline{c} ,\overline{d} ) = (c\pmod{p^n}, d\pmod{p^n}).$$ %
L'application de ${\Gamma }_0 (p^n )\backslash SL_2(\grdZ)$ vers $H_1 ({X_0}
({p^n}),\ {\rm pointes} ; \grdZ ) $ : $ g \mapsto \{ g\cdot 0,g\cdot
\infty \} $ s'identifie alors a une application de ${\grdP}^1 \left( \grdZ /
{{p^n}\grdZ} \right)$ vers la m\^{e}me chose, qu'on note $\xi $ :
$$\xi (\overline{w} ,\overline{t} ) = \left\{ \overline{ \left(
\begin{array}{cc}
a & b \\
w & t
\end{array}
\right) \cdot 0 } , \overline{\left(
\begin{array}{cc}
a & b \\
w & t
\end{array}
\right) \cdot \infty  } \right\} = \left\{  \frac{b}{t} , \frac{a}{w}
\right\}  , $$
avec $w,t$, rel\`{e}vements dans $\grdZ$ de $\overline{w} $ et $\overline{t}
\,  \in  \grdZ / {p^n \grdZ} $ et $a,b\, \in \grdZ $ tels que $\left(
\begin{array}{cc}
a & b \\
w & t
\end{array}
\right)$ soit dans $SL_2 (\grdZ)$.

   On sait de plus qu'on a :
 $$T_r \{ 0,\infty \} = \sum_{\begin{array}{c}
0 \leq  w <  t \\
0 \leq  v <  u \\
ut-vw  =r
\end{array} }
\xi (\overline{w} ,\overline{t} )\ ,$$ %
o\`{u} on pose $\xi (\overline{w} ,\overline{t} ) = 0$ si pgcd$(w,t,p) >1$
(voir \cite{Artin}, th\'eor\`eme 2 et proposition 20, ou la d\'emonstration du
lemme 2 de \cite{merel}).

  Soit $\sigma = \overline{\left( \begin{array}{rc}
0 & 1 \\
-1 & 0
\end{array}
\right) } $ et $\tau = \overline{\left( \begin{array}{cc}
0 & -1 \\
1 & -1
\end{array}
\right) }$. On choisit pour repr\'{e}sentants de ${\grdP }^1 (\grdZ / {p^n}
\grdZ )$ : $\{ (R_1 ,1),\ R_1 $ un syst\`{e}me de repr\'{e}sentants de
$\grdZ /{p^n}\grdZ \}\ \cup \{ (1,p.R_2 ),\ R_2 $ un syst\`{e}me de
repr\'{e}sentants de $\grdZ /{p^{n-1}} \grdZ \}$. On note $w/t $ au lieu de
$(\overline{w} ,\overline{t} )$, souvent. On fait agir $SL_2 (\grdZ )$ sur
${\grdP }^1 (\grdZ /{p^n} \grdZ )$ \`{a} droite, comme d'habitude. En
particulier :
$$(\overline{w} ,\overline{t} ) \cdot \sigma = (\overline{w} ,\overline{t} )
\overline{\left( \begin{array}{cc}
0 & -1 \\
1 & 0
\end{array}
\right) }  =(\overline{-t} ,\overline{ w} ),$$ %
et de m\^{e}me $(w / t )\cdot \tau = -t / (w+t ) $.

  On note encore ${\Sigma }_r  = \{ (\overline{w} ,\overline{t} ),\ 0\leq w<t
\ /\ $il existe $(v,u),\ 0\leq v<u\ ,\ 0\leq (ut-vw)\leq r\} \backslash \{
(\overline{1} , \overline{r} )\} $, et $\grdZ [ {\grdP }^1 (\grdZ /{p^n}
\grdZ ) { ] }^{\sigma }$ (respectivement, $\grdZ  [ {\grdP }^1 (\grdZ /{p^n}
\grdZ ) {] }^{\tau } $) d\'{e}signe l'ensemble des \'{e}l\'{e}ments de
$\grdZ [ {\grdP }^1 (\grdZ /{p^n} \grdZ  ) ]$ stables par l'action de $\sigma$
(respectivement, $\tau $). Le symbole $\sum_{\sigma }^{ } $ d\'{e}signera une
somme \`{a} valeurs dans le $\grdZ$-module
$\grdZ [ {\grdP }^1 (\grdZ /{p^n} \grdZ )
{ ] }^{\sigma }$, et de m\^{e}me avec~$\tau $.

  On a la suite exacte (voir par exemple \cite{Mer}) :
$$\grdZ [ {\grdP}^1 \left( \grdZ / {p^n}\grdZ \right)  { ] }^{\sigma } \times
\grdZ [{\grdP}^1 \left( \grdZ / {p^n}\grdZ \right) { ] }^{\tau }
\stackrel{{\phi}_1}{\rightarrow} \grdZ [ {\grdP}^1  (\grdZ / {p^n} \grdZ ) ]
\stackrel{{\phi}_2}{\rightarrow} H_1 (X_0 ({p^n}),\, {\rm pointes}\, ;\,
\grdZ ) \to 0,$$
avec ${\phi}_1$ : $(\Sigma {\alpha}_x x,\Sigma {\beta }_x x ) \mapsto
(\Sigma {\alpha }_x x + \Sigma {\beta }_x x),$ et ${\phi}_2$ : $\Sigma
{\lambda}_x x \mapsto \Sigma {\lambda}_x .\xi (x)$.
\subsection{Preuve de la proposition.}
Posons $D=sd$. Supposons $p^n >C^2 .D^6$, et que pour un $r\leq D$, on ait
une relation de liaison : $\sum_{i=1}^{r} {\lambda }_i {T_i} \{ 0,\infty \}
=0$ dans $H_1 (X_0 ({p^n} ),\, {\rm pointes}\, ;\, {\grdF }_l )$. On va
montrer que ${\lambda }_r$ est nul : ce qui suffira \`a la preuve de la
proposition. Ce qui pr\'{e}c\`{e}de permet d'\'{e}crire :
$$\sum_{i=1}^{r} {\lambda }_i T_i \{ 0,\infty \} =0 \iff {\lambda }_r
.\left( 1/r \right) +
\sum_{{\Sigma }_r} {\mu }_{(\overline{w} ,\overline{t} )} \left( w/t
\right) = \sum_{\sigma } {\alpha }_x (x) -\sum_{\tau } {\beta }_{x} (x).$$
De m\^{e}me que dans \cite{merel} on prouve avec l'aide de Fouvry le :
\begin{lemm}
Soit $A$ et $B$ deux intervalles de $\{ 1,2,...,{p^n}-1 \}$ tel que
$$|A|.|B| \geq C'.{p^{3n/2}} ,$$
o\`u $C'=8$ si $p$ est impair, et $C'=8\sqrt{2}$ si $p=2$. Alors il existe
$y\in A$ et $z\in B$ tel que $y.z=\ -1\pmod{p^n}$.
\end{lemm}
(On a prouv\'e une version moins bonne de ce lemme, utilisant une constante
$C=(512{\pi}^2 )/(2\sqrt{2} -1)$ ; la forme ici utilis\'ee de ce lemme est
obtenue en optimisant les calculs par Oesterl\'e dans \cite{Oesterle}.)
On consid\`{e}re le graphe de ${\grdP }^1 (\grdZ /{p^n}\grdZ )$ dont les
ar\^{e}tes sont l'action de $\sigma $ et $\tau $ (voir {\bf Figure 1}) : on a
$(\tau \sigma )= \overline{\left( \begin{array}{cc}
-1 & 0 \\
-1 & -1
\end{array}
\right)}  ,$ donc $(\overline{w} ,\overline{t} )\cdot \tau \sigma = (
\overline{w} ,\overline{t} ) +\overline{1}$ (et de m\^eme $(\overline{w} ,
\overline{t} )\cdot \sigma {\tau }^2 =(\overline{w} ,\overline{t} )-
\overline{1} ).$

\setlength{\unitlength}{0.7cm}

\begin{picture}(26,20)(3,0)
\thicklines
\put (10,16){\framebox (2,2){$(1,r)$}}
\thinlines
\put (9.7,17){\vector (-1,-1){1.8}}
\put (8,16){${\tau }^2$}
\put (5,14){\framebox (3,1){$(-r-1,1)$}}
\put (7.3,13.7){\vector (1,-1){1}}
\put (8,13.2){$\sigma$}
\put (8,11){\framebox (2.5,1){$(1,r+1)$}}
\put (7,10.5){${\tau}^2$}
\put (7.7,10.7){\vector (-1,-1){1}}
\put(6.5,9){.}
\put(7,8.5){.}
\put(7.5,8){.}
\put(7.9,7.5){.}
\put(8,6.5){\vector (-1,-1){1}}
\put (6.8,6.2){${\tau}^2$}
\put (5,4){\framebox (2,1){($y,1)$}}
\put(10,5.5){$\sigma$}
\thicklines
\put (8,4.5){\vector (4,1){4}}
\thinlines
\put (13.2,16.8){$\sigma$}
\put (15,15){\framebox (2,1){$(-r,1)$}}
\put (12.6,17){\vector (2,-1){2}}
\put (17,14.7){\vector (-1,-2){1}}
\put (16,11){\framebox(2,1) {...}}
\put (18.3,10.7){\vector (1,-2){0.9}}
\put (19,8.2){.}
\put(19,7.7){.}
\put(19,7.2){.}
\put(18,7){\vector(-3,-1){2.7}}
\put(13,5){\framebox(2,1){$(z,1)$}}
\put (4,15){\vector(0,-1) {10.5}}
\put(2.5,9){${\cal A}$}
\put(4,3.5){.}
\put(4,3){.}
\put(4,2.5){.}
\put(20.5,16){\vector (0,-1){11}}
\put(21,10){${\cal B}$}
\put(20.5,4.5){.}
\put(20.5,4){.}
\put(20.5,3.5){.}
\put(11.5,2){\bf Figure 1}
\end{picture}
On va montrer qu'on a sur ce graphe, ``de part et d'autre'' de $(\overline{1}
,\overline{r} )$ (c'est-\`{a}-dire, contenant $(\overline{1} ,\overline{r} )
\cdot {\tau }^2 =
(\overline{-r-1} ,\overline{1} )$ et $(\overline{1} ,\overline{r} )\cdot
\sigma =(\overline{-r} ,\overline{1} )$ respectivement), deux chemins
${\cal A}$ et ${\cal B}$ ne rencontrant pas d'\'{e}l\'{e}ments de
${\Sigma }_r$, (les autres \'{e}l\'{e}ments du graphe intervenants dans $\sum
{\lambda }_i T_i \{ 0,\infty \} )$, et qui contiennent des intervalles de
cardinal sup\'{e}rieur \`{a} $({{p^n}/D})-D-2$ et $({{p^n}/{D^2}})-2$
respectivement. Alors par le lemme de th\'eorie analytique des nombres, pour
$$(({{p^n}/D})-D-2).(({{p^n}/{D^2}})-2)\geq C'.{p^{3n/2}} ,\ i.e.$$
$${p^n}\geq C^2 . {D^6} ,$$
on aura dans ${\cal A}$ et dans ${\cal B}$ respectivement des \'el\'ements $y$
et $z$ tel que $y\cdot \sigma =\frac{-1 }{y} =z$ (on explicite ces calculs,
et notament le passage de $C'$ \`a $C$, \`a la fin de la section). De plus
$y\cdot \sigma $ sera un \'{e}l\'{e}ment de ${\cal A}$, puisque ce chemin est
de forme :
$$\begin{array}{rcl}
...\to (\overline{a} ,\overline{1} ) \stackrel{\tau }{\rightarrow } & (
\overline{-1} ,\overline{a+1} ) &
\stackrel{\sigma }{\rightarrow} (\overline{a+1} ,\overline{1} ) \to... \\
\longrightarrow  & +\overline{1} =\tau \sigma & \longrightarrow
\end{array}$$
Les deux chemins se ``rencontrent'', donc. Or on a, pour tout \'el\'ement $x$
du graphe, ${\mu }_x = {\alpha }_x - {\beta }_x$ par ce qui pr\'{e}c\`{e}de.
Puisque les deux chemins ne rencontrent pas $\Sigma _r$, si $x$ en est, on a
${\mu }_x =0$ donc ${\alpha }_x ={\beta }_x$. De plus, on parcourt ces chemins
en faisant agir $\sigma $ ou $\tau $ ; donc si $x'=x\cdot \sigma $,
${\alpha}_{x'}={\alpha }_{x\cdot \sigma }={\alpha }_x ={\beta }_x =
{\beta }_{x'}$, de m\^{e}me avec $\tau $ - donc ${\alpha }_x ={\beta }_x
\equiv {\alpha }_{-r} $ sur les deux chemins. Mais
$${\lambda }_{\overline{r} }=
{\mu }_{\frac{1}{r} } = {\alpha }_{\frac{1}{r} } -{\beta }_{\frac{1}{r} } =
{\alpha }_{\frac{1}{r} \cdot \sigma } - {\beta }_{\frac{1}{r} \cdot {\tau }^2}
=0.$$
Montrons donc l'existence de ces chemins. (Dans tous les calculs qui
suivent, on confond l'\'ecriture d'un entier $w$ et de sa r\'eduction
$\overline{w}$, pour all\'{e}ger les notations).

$${ }$$
\underline{Premier chemin : ${\cal A}$} partant de
$\frac{1}{r} \  {\tau }^2 =-r-1.$  \\

1) Si $\frac{w}{t} -(-r-1)= \overline{a}$, avec $a$ choisi dans
$\{ -{p^n}+1,...,-1,0 \}$, et $\overline{a} $ : classe de $a$ mod $\ {p^n}$
($\iff w+t(r+1)=at+b.{p^n}$, $b\in \grdZ $).\\

\underline{Si $b=0$} : $\ t(r+1)-at=-w $ : incompatibilit\'{e} de signes. \\

\underline{Si $b \neq 0$} : $\ |a|=\frac{1}{t} |b.{p^n}-t(r+1)-w|\geq
\frac{({p^n}-D(D+1)-D)}{D} \geq \frac{p^n}{D} -D-2 $
(on a en effet : det$\left( \begin{array}{cc}
u & v \\
w & t
\end{array}
\right)=k\leq r,\ u>v\geq 0,\ t>w\geq 0$, donc :
$k=ut-vw\geq ut-(u-1)(t-1)$ et
$u+t-1\leq r,\ t\leq r-u+1 \leq r \leq D$).\\

2) Si $(\frac{w}{t} )\sigma -(-r-1)= \frac{-t}{w} +r+1=\overline{a} ,\ a\in
\{-{p^n}+1,...,-1,0\}\ $; $\ -t+w(r+1)=aw+b{p^n}\ $; \\

\underline{Si $b=0$} : $\ w(r+1-a)=t\ ;$ mais $(r+1-a)\geq r+1,$ et $0\leq w<t
\leq r$ : contradiction.\\

\underline{ Si $b\neq 0$}, $\ |a|\geq ( \frac{{p^n}-D(D+1)+D}{D} )\geq
\frac{p^n}{D} -D$.\\

 On peut donc ``reculer'' (...$\ \alpha \stackrel{\sigma }{\to } .
\stackrel{\tau ^2}{\to } \alpha -1 \stackrel{\sigma }{\to } .
\stackrel{\tau ^2}{\to } \alpha -2\ ...$) \`{a} partir de $(-r-1)$, et
d\'{e}crire ainsi un chemin contenant un intervalle de cardinal sup\'{e}rieur
\`{a} ${p^n}/{D}\ -D-2$.\\

\underline{Second chemin.} On doit l\`{a} distinguer deux cas, selon que $p$
divise ou non $r$ (voir {\bf Figure 2}).\\

{\bf a}) Si $p$ ne divise pas $r$, chemin ${\cal B}$ : on part de
$(\frac{1}{r} )$ lui m\^{e}me, on recule de m\^{e}me : \\

1) $\frac{w}{t} -\frac{1}{r} =\overline{a} ,\ -{p^n} < a\leq 0 \iff wr-t=art+
b.{p^n}\  $;\\

\underline{$b=0$}\ : $t=r(w-at) \Rightarrow a=0,\ w=1,\ t=r$ : c'est
$\frac{1}{r} $ lui-m\^{e}me.\\

\underline{$b\neq 0$}\ : $|a|\geq \frac{{p^n}-{D^2}}{D^2}$ de m\^{e}me que
plus haut.\\

2) $-\frac{t}{w} -\frac{1}{r} =\overline{a} $ : $-rt-w=awr+b.{p^n}$\ ; \\

\underline{$b=0$} $\Rightarrow r|w,$ impossible (car $w\leq r-1,$ et : $w=0
\Rightarrow t=0$).\\

\underline{$b\neq 0$} $\Rightarrow |a|\geq \frac{p^n}{D^2} -2$.\\

{\bf b}) Si $p$ divise $r$, on a alors que le chemin ${\cal B}$
pr\'{e}c\'{e}dent : $\frac{1}{r} \stackrel{\sigma {\tau }^2}{\longrightarrow}
\frac{1-r}{r} \stackrel{\sigma \tau ^2}{\longrightarrow } ...$ est bien de
longueur sup\'{e}rieure \`{a} $(\frac{p^n}{D^2} )$ ; mais il ne contient pas
cette fois d'intervalle, puisque $r$ n'est pas inversible modulo ${p^n}$, donc
l'action de $\sigma \tau ^2 $ ne correspond plus \`{a} l'addition de $(-1)$
(les $k/pl$ ne sont plus relevables en \'{e}l\'{e}ments de $\grdZ /{p^n}
\grdZ$). Cependant, le calcul pr\'{e}c\'{e}dent montre que $\frac{r}{r-1} =
\frac{1}{r} \cdot \sigma {\tau ^2}\sigma $\ ; et $(r-1)$ est, cette fois,
inversible modulo $p^n$, donc en ``avan\c{c}ant''  \`{a} partir de cet
\'{e}l\'{e}ment et \`{a} l'aide de $(\tau \sigma )$, on aura bien un
intervalle ; on note ${\cal B}'$ ce chemin. On minore encore une fois sa
longueur : \\

1) $\frac{w}{t} -\frac{r}{r-1}= \overline{a} \iff w(r-1) -rt= at(r-1)+b.{p^n}$
 (on \'{e}crit cette fois cela avec $0\leq a\leq {p^n}-1$).\\

\underline{Si $b=0$} :\  $\ t\left( r+a(r-1) \right) =w(r-1)$ ; mais $w<t,$
et $a(r-1)+r>r-1$, contradiction. \\

\underline{Si $b\neq 0$} : $\ |a|\geq \frac{p^n}{D^2} -1.$\\

2) $\frac{w}{t} \cdot \sigma -\frac{r}{r-1} =-\frac{t}{w} -\frac{r}{r-1} =
\overline{a} \iff -t(r-1)-rw=aw(r-1) + b.{p^n}$. \\

\underline{$b=0$} : $\ t(r-1)=-w \left( a(r-1)+r \right) $, contradiction de
signes. \\

\underline{ $b\neq 0$} : $\ |a| \geq \frac{p^n}{D^2} -2.$ \\

 Dans chaque cas, on obtient bien que ce second chemin contient un intervalle
de cardinal sup\'{e}rieur ou \'{e}gal \`{a} $({p^n}/{D^2}) -2$.$\Box$

\begin{picture}(20,20.5)(1.5,-2.5)

\thinlines
\put (3,16.3){\vector (0,1){1}}
\put (3.2,16.5){${\tau }^2$}
\put(3,17.5){.}
\put(3,17.7){.}
\put(3,17.9){.}

\thicklines

\put (2,14){\framebox (2,2){$(1,r)$}}
\thinlines
\put (3.6,13.4){$\sigma$}
\put (4,13.8){\vector (1,-1){0.7}}

\put (5,12){\framebox (2,1){$(-r,1)$}}
\put (5,11.8){\vector (-1,-1){0.7}}
\put (3.8,11.5){${\tau}^2$}
\thicklines
\put (2,10){\framebox (2.5,1){$(1-r,r)$}}
\thinlines
\put (4,9.8){\vector (1,-1){0.7}}
\put (4.5,9.5){$\sigma$}
\put (5,8){\framebox (2.4,1){$(r,r-1)$}}

\put (5,7.8){\vector (-1,-1){0.7}}
\put (4.3,7.5){${\tau}^2$}
\put (1.5,6){\framebox (2.7,1){$(1-2r,r)$}}
\put (4,5.8){\vector (1,-1){0.7}}
\put (4,5){$\sigma$}
\put (5.1,4){\framebox (2.7,1){$(r,2r-1)$}}
\put (5,3.8){\vector (-1,-1){0.7}}
\put (5,3){${\tau}^2$}
\put (1.6,2){\framebox (2.7,1){$(1-3r,r)$}}
\put (4,1.8){\vector (1,-1){0.7}}
\put(5,0.5){. . . }

\put (7.5,8.4){\vector (1,0){1}}
\put (7.8,8.5){$\tau$}
\put (8.9,8){\framebox (3.6,1){$(1-r,2r-1)$}}
\put (10.4,7.8){\vector (0,-1){0.7}}
\put (10.6,7.5){$\sigma$}
\put (8.9,6){\framebox (3.6,1){$(2r-1,r-1)$}}
\put (12.7,6.4){\vector (1,0){1}}
\put (13,6.7){$\tau$}
\put (13.9,6){\framebox (3.6,1){$(1-r,3r-2)$}}
\put (15.3,5.8){\vector (0,-1){0.7}}
\put (15.5,5.3){$\sigma$}
\put (13.9,4){\framebox (3.6,1){$(3r-2,r-1)$}}
\put (17.7,4.5){\vector (1,0){1}}
\put (18,4.7){$\tau$}
\put (19,4.3){. . . .}

\put(1,14){\vector (0,-1){12.5}}
\put(0,8){${\cal B}$}

\put(7,7){\vector (1,-1){6}}
\put(9,4){${\cal B}'$}

\put(8,-1){\bf Figure 2}

\end{picture}
On montre pour finir comment on passe du lemme 5.3 \`{a} la proposition
1.7 : on a vu qu'on pouvait prendre $|A|\geq ({p^n}/{D^2} )-2$ et $|B|\geq
({p^n}/D) -D-2 $. On a suppos\'e qu'\'etait satisfaite la condition de la
proposition : $p^n \geq (C^2).D^6$. On a $C^2\geq 65$, $s\geq 2$, et les
bornes du corollaire 1.8 sont sup\'erieures (!) \`a celles qu'on connaissait
d\'ej\`a pour les degr\'es 1 et 2, donc on peut supposer $D\geq 6$. Minorons
la taille de $A$ : on a $p^n /D^2 \geq 65. D^4 \geq 84240=42120.(2)$, donc
$|A|\geq (42119/42120).(p^n /D^2 )$. Pour $B$, on a : $D+2\leq (4/3).D$, et
$p^n /D\geq 65.D^5 .D\geq 379080.(D+2)$ ; donc $|B|\geq (379079/379080).(p^n
/D)$. Si on pose $\lambda :=(42119/42120).(379079/379080)$, on a donc :
$$|A|.|B|\geq \lambda .p^{2n} /D^3 .$$
Pour que les conditions du lemme soient v\'erifi\'ees, il suffit qu'on ait
$\lambda .p^{2n} /D^3 \geq C'.p^{3n/2}$, {\em i.e.} $p^n \geq (C'^2 /
{\lambda}^2 ).D^6$, et \'evidemment $C$ a \'et\'e choisie pour que \c{c}a
marche. $\Box$ %
$${ }$$
{\bf Remerciements.}\  Je remercie ici Lo\"{\i}c Merel qui m'a sugg\'{e}r\'{e}
l'id\'{e}e essentielle de ce papier, et Bas{ }Edixhoven qui encadre tout ce
travail. Merci aussi \`{a} Joseph Oesterl\'{e} pour l'id\'ee du lemme 1.4, et
\`{a} \'Etienne Fouvry pour le lemme de th\'{e}orie analytique des nombres.\\

\begin{picture}(10,6)
\put(0,5){\line(1,0){6}}
\put(0,4){Pierre Parent}
\put(0,3){IRMAR}
\put(0,2){Universit\'e de Rennes}
\put(0,1){35 042 RENNES C\'edex France}
\put(0,0){E-mail : {\tt parent@@clipper.ens.fr}}
\end{picture}


\begin{thebibliography}{XXXXXXXXXXXXXXXXXX56}

\bibitem[Atkin-Lehner 70]{Atkin-Lehner}

{A.O.L. Atkin, J. Lehner, {\it Hecke Operators on ${\Gamma}_0 (m)$}, Math.
Ann. 185, pages 134-160, (1970).}


\bibitem[Bump-Friedberg-Hoffstein 90]{Bump-Friedberg-Hoffstein}

{D. Bump, S. Friedberg, J. Hoffstein, {\em Nonvanishing theorems for
L-functions
of modular forms and their derivatives}, Inventiones Mathematic\ae  102, pages
543-618, (1990).}


\bibitem[Coleman-Edixhoven 96]{Coleman-Edixhoven}

{Robert F. Coleman, Bas Edixhoven, {\em On the semi-simplicity of the
$U_p$-operator on modular forms}, article en pr\'eparation.}


\bibitem[Diamond-Im]{Diamond}

{Fred Diamond, John Im, {\it Modular forms and modular curves}, Canadian
Mathema\-tical Society Conference Proceedings, \`a para\^{\i}tre.}

\bibitem[Edixhoven 94]{Bas}

{B. Edixhoven, {\it Rationnal torsion points on elliptic curves over number
fields}, As\-t\'{e}\-ris\-que (S\'{e}minaire Bourbaki 782) (1993).}

\bibitem[EGA IV]{egaiv}

{A. Grothendieck, {\it \'El\'ements de g\'eo\-m\'e\-trie alg\'ebrique IV,
(\'Etude locale des morphismes de sch\'emas)}, Publications math\'ematiques de
l'I.H.E.S. 28 (1965).}

\bibitem[Hartshorne]{Hartshorne}

{R. Hartshorne, {\it Algebraic Geometry}, GTM 52, Springer-Verlag.}

\bibitem[Hooley]{hool}

{Hooley, {\it Applications of Sieve methods to Number Theory}, Cambridge
University Press, page 35 (1935).}


\bibitem[Kamienny 92a]{kamienny inv}

{S. Kamienny, {\it Torsion points on elliptic curves and $q$-coefficients of
modular forms}, Inventiones Mathematic\ae  109, pages 221-229 (1992).}

\bibitem[Kamienny 92b]{kamienny}

{S.\ Kamienny, {\it Torsion points on el\-lip\-tic cur\-ves over fields of
higher degree}, In\-ter\-na\-tio\-nal Mathe\-ma\-tics Re\-search No\-ti\-ces
6 (1992).}


\bibitem[Kamienny-Mazur]{kamienny-mazur}

{S. Kamienny et B. Mazur, {\it Rationnal torsion of prime order in elliptic
cur\-ves over number fields}, Ast\'{e}risque (\`{a} para\^{\i}tre).}

\bibitem[Katz-Mazur]{Katz-Mazur}

{N. Katz and B. Mazur, {\it Arithmetic Moduli of Elliptic Curves}, Annals of
Mathematics Studies 108, Princeton University Press (1985).}

\bibitem[Kolyvagin-Logachev 90]{Kolyvagin}

{V.A. Kolyvagin, D.Yu. Logachev, {\it Finiteness of the Shafarevich-Tate group
and the group of rationnal points for some modular abelian varieties},
Leningrad Math. J., Vol. 1, number 5 (1990).}

\bibitem[Lang]{Lang}
{Serge Lang, {\it Introduction to modular forms}, Grundl. Math. Wiss. 222,
Springer-Verlag (1976).}

\bibitem[Manin 69]{man}
{Manin, {\it A Uniform Bound for p-torsion of Elliptic Curves}, Izv. Akad.
Nau. CCCP 33 (1969).}


\bibitem[Manin 72]{manin1}

{Y. Manin, {\it Parabolic points and zeta function of modular curves}, Math.
USSR Izvestija 6, pages 19-64 (1972).}

\bibitem[Mazur 77]{mazur}

{B. Mazur, {\it Modular curves and the Ei\-sen\-stein ideal}, Publications
math\'{e}matiques de l'I.H.E.S. 47, pages 33-186 (1977).}

\bibitem[Mazur 78]{rational}

{B. Mazur, {\it Rational Isogenies of Prime De\-gree}, Inventiones
Mathematicae
44, pages 129-162 (1978).}

\bibitem[Merel 93]{Mer}

{Lo\"{\i}c Merel, {\it Sur quelques aspects g\'{e}o\-m\'{e}\-tri\-ques et
arithm\'{e}tiques de la th\'{e}orie des symboles modulaires}, Th\`{e}se de
doctorat de l'universit\'{e} de Paris VII  (1993).}

\bibitem[Merel 94]{Artin}

{Lo\"{\i}c Merel, {\em Universal Fourier expansions of modular forms}, in {\em
On Artin's conjecture for odd 2-dimensional representations}, Lecture Notes in
Mathematics 1585, Springer-Verlag (1994).}

\bibitem[Merel 95]{merel}

{Lo\"{\i}c Merel, {\it Bornes pour la torsion des courbes elliptiques sur les
corps de nombres}, Inventiones Mathematic\ae , \`{a} para\^{\i}tre.}

\bibitem[Mumford]{mumford}

{David Mumford, {\it Abelian Varieties}, Oxford University Press (1970).}

\bibitem[Murty-Murty 91]{Murty-Murty}

{M. R. Murty, V. K. Murty, {\em Mean values of derivatives of modular
L-series,}
Ann. Math. 133, pages 447-475, (1991).}

\bibitem[Oesterl\'{e} 94]{Oest}

{Joseph Oesterl\'{e}, {\it Torsion des courbes elliptiques sur les corps de
nombres}, article \`{a} para\^{\i}tre.}

\bibitem[Oesterl\'e 96]{Oesterle}

{J. Oesterl\'e, {\it Un lemme de th\'eorie analytique des nombres}, \`a
para\^{\i}tre.}

\bibitem[Sali\'{e} 31]{salie}

{Sali\'{e}, {\it \"{U}ber die Kloostermanschen Summen}, Math. Zeit. 34, pages
91-109 (1931).}

\bibitem[SGA 7-I]{SGA7I}

{Grothendieck et al., {\it S\'eminaire de g\'eom\'etrie alg\'ebrique du
Bois-Marie 7-I (Grou\-pes de mo\-no\-dro\-mie en g\'eo\-m\'e\-trie
al\-g\'e\-bri\-que)},
Lec\-tu\-re No\-tes in Mathe\-ma\-tics 288, Springer-Verlag, (1972).}

\bibitem[Serre]{serre}

{Jean-Pierre Serre, {\it Groupes alg\'ebriques et corps de classes,} Hermann
(1959).}

\bibitem[Shimura]{Shimura}

{Goro Shimura, {\it Introduction to the Arith\-metic Theory of Automorphic
Functions}, Princeton University Press (1971).}

\end{thebibliography}
\end{document}